\documentclass[11pt]{article}
\usepackage{amssymb,amsfonts,amsmath,amsthm,cite}
\usepackage{epsfig}
\newtheorem{thm}{Theorem}[section]

\newtheorem{lem}[thm]{Lemma}

\newtheorem{exam}[thm]{Example}

\makeatletter \@addtoreset{equation}{section}

\def\qed{\hfill \rule{4pt}{7pt}}
\def\pf{\noindent {\it Proof.}\quad}

\title{\bf\Large  Stanley's Lemma and  Multiple Theta Functions }
\author{William~Y.~C.~Chen$^{a,b}$ and Lisa~H.~Sun$^a$
\date{$^{a}$Center for Combinatorics, LPMC\\
 Nankai University, Tianjin 300071, P. R. China\\[5pt]
$^{b}$Center for Applied Mathematics \\
Tianjin University, Tianjin 300072, P. R. China
\vskip 0.2 cm
emails: chen@nankai.edu.cn, sunhui@nankai.edu.cn}}

\begin{document}
\maketitle

\allowdisplaybreaks
\noindent {\bf Abstract.} We present an algorithmic approach to
 the verification of identities on multiple theta functions in the form of products of theta functions $[(-1)^{\delta}a_1^{\alpha_1}a_2^{\alpha_2}\cdots a_r^{\alpha_r}q^{s}; q^{t}]_\infty$, where $\alpha_i$ are integers, $\delta=0$ or $1$, $s\in  \mathbb{Q}$, $t\in \mathbb{Q}^{+}$, and the exponent vectors $(\alpha_1,\alpha_2,\ldots,\alpha_r)$  are linearly independent over $\mathbb{Q}$.
For an identity on such multiple theta functions, we provide an algorithmic approach
 for computing a system of contiguous relations satisfied by all the involved multiple theta functions. Using Stanley's Lemma on the fundamental parallelepiped,  we show that
 a multiple theta function can be determined by a finite number of its coefficients. Thus such an identity can be reduced to a finite number of simpler relations. Many classical multiple theta function identities fall into this framework, including Riemann's addition formula and the extended Riemann identity.



\section{Introduction}
\allowdisplaybreaks

Theta functions arise in the  study of  the Riemann zeta functions, the Weierstrass elliptic functions and the $q$-gamma functions. There are many identities on the classical Jacobi theta functions including Jacobi's triple product identity, the quintuple product identity and the septuple product identity due to Farkars and Kra \cite{far08, farkra99}. For more properties of theta functions, one can see \cite{adiberbhawat85, AAR99, berndt06, whiwat27}.

 Assume that $|q|<1$.
The theta function $\theta(z)$ is defined by
\begin{equation}\label{thetaf}
\theta(z)=[z;q]_\infty=(z,q/z;q)_\infty,
\end{equation} where the $q$-shifted factorial $(z;q)_\infty$ is given by
\begin{equation*}
(z;q)_\infty=\prod_{k=0}^\infty (1-zq^k).
\end{equation*}
A multiple theta function is defined by
\begin{equation}\label{notation}
[z_1, z_2,\ldots, z_n;q]_\infty=[
z_1;q]_\infty[z_2;q]_\infty \cdots[z_n;q]_\infty.
\end{equation}

In this paper, we are concerned with identities on a special class of multiple theta functions.
Assume that  $a_1, a_2, \ldots, a_r$ are complex variables.
For a vector $\alpha=(\alpha_1, \alpha_2, \ldots, \alpha_r)$
of integers, we adopt the common notation $a^\alpha=a_1^{\alpha_1}a_2^{\alpha_2}\cdots
a_r^{\alpha_r}$.
We shall restrict our attention to multiple theta functions  of the form:
\begin{equation}\label{thetaeq}
\theta (a)=a^{\tau} \prod\limits_{i=1}^m f_i(a),
\end{equation}
where $\tau \in \mathbb{Q}^r$,  $1\le m\le r$ and
\begin{equation}\label{factormt}
f_i(a)=[(-1)^{\delta_i}a^{\gamma_{i}}q^{s_i}; q^{t_i}]_\infty,
\end{equation}
such that for $1\leq i \leq m$,  $\gamma_{i}=(\gamma_{i,1},\gamma_{i,2},\ldots,\gamma_{i,r})$ is
 a vector
of integers, $\delta_i =0$ or $1$, $s_i \in  \mathbb{Q}$ and $t_i\in \mathbb{Q}^{+}$.
 We further assume that the exponent vectors $\gamma_1, \gamma_2,\ldots,\gamma_m$  are linearly independent and $\tau$ is a linear combination of $\gamma_1, \gamma_2,\ldots,\gamma_m$ over $\mathbb{Q}$.  For example,
\begin{equation}\label{examtheta}
\theta(a,b,c,d)=(ab,q/ab;q)_\infty(a/b,bq/a;q)_\infty (cd,q/cd;q)_\infty (c/d,qd/c;q)_\infty
\end{equation}
is a multiple theta function in the form of \eqref{thetaeq}. The
 exponent vectors in the factors of $\theta(a,b,c,d)$ are
  \[
  \gamma_1=(1,1,0,0),\ \gamma_2=(1,-1,0,0),\ \gamma_3=(0,0,1,1),\ \gamma_4=(0,0,1,-1),
  \]
  which are  linearly independent over $\mathbb{Q}$.

 Many identities on the multiple theta functions arise as generalizations of theta function identities. Winquist \cite{Win69} found an  identity on  bivariate theta functions which plays an important role in proving Ramanujan's congruence for the partition function modulo 11.   Carlitz and Subbarao \cite{CarSub72}, and Hirschhorn \cite{Hir86} obtained further generalizations of Winquist's identity. Liu \cite{Liu07} derived  an addition formula for the Jacobi theta functions  by using the theory of elliptic functions, which specializes  to the Ramanujan cubic theta function identity and Winquist's identity.  Ewell \cite{Ewell95}
 found a sixfold infinite-product identity on multiple theta functions with three variables. Shen \cite{shen94} obtained a collection of addition formulae of theta functions by using Fay's trisecant identity, which plays a vital role in the study of Riemann surfaces (see, for example \cite{farkra01}). Bailey \cite{Bai36} deduced an identity on four multiple theta functions with five variables by applying the basic hypergeometric series $_8\phi_7$. Slater \cite{Sla54} extended this identity  to seven multiple theta functions with six variables by using $_{10}\phi_9$ series.
Malekar and Bhate \cite{MalBha12} employed  the discrete Fourier transform to derive several fourth order identities on the Jacobi theta functions such as the extended
Riemann identity. Recently, Cao \cite{cao11} established a correspondence between identities on the Jacobi theta functions and integer matrix exact covering systems, which can be used to produce identities on products of Ramanujan's theta functions.

For a  multiple theta function $\theta (a_1,a_2,\ldots,a_r)$, a contiguous relation is meant to be a relation of the form
\begin{align}\label{recmulti}
\frac{\theta(a_1q^{x_1},a_2q^{x_2},\ldots,a_rq^{x_r})}{\theta(a_1,a_2,\ldots,a_r)}
=\frac{(-1)^{\rho}}{a_1^{v_{1}}a_2^{v_{2}}\cdots
a_r^{v_{r}}q^{s}},
\end{align}
where $\rho=0$ or 1, $x_{i}, s\in \mathbb{Q}$ and $v_{i} \in \mathbb{Z}$.
 The vector $v=(v_{1},v_{2},\ldots,v_{r})$ in the denominator  of the contiguous relation \eqref{recmulti}  is referred to as the exponent vector.
Note that once the exponent vector $(v_{1},v_{2},\ldots,v_{r})$ is determined, there may be different vectors $(x_1,x_2,\ldots,x_r)$ satisfying the contiguous relation (\ref{recmulti}). Nevertheless, if we wish to derive
a recurrence relation on the  coefficients of $\theta(a)$  from
  the contiguous relation (\ref{recmulti}), we mainly need the exponent vector
  $(v_1,v_2,\ldots,v_r)$.
 For example, for $\theta(a,b,c)=(a,q/a;q)_\infty(a/bc,bcq/a;q)_\infty$, we have
 \begin{equation}\label{dcont}
 \frac{\theta(aq,bq,c)}{\theta(a,b,c)}=\frac{\theta(aq,b,cq)}{\theta(a,b,c)}
 =-\frac{1}{a}.
 \end{equation}
Note that $\theta(a,b,c)$ can be expanded by using Jacobi's triple product identity
\begin{align}\label{jactri}
(q,z,q/z;q)_\infty=\sum_{k=-\infty}^\infty (-1)^k q^{k \choose 2} z^k.
\end{align}
 Let
 \[
 \theta(a,b,c)=\sum_{(n,m,k)\in \mathbb{Z}^3} h(n,m,k)a^nb^mc^k.
  \]
Applying \eqref{jactri}, we see that the above sum ranges over integer vectors
$(n,m,k)$ such that
\[
  (n,m,k)=\ell_1 (1,0,0)+\ell_2 (1,-1,-1),
\]
where  $\ell_1, \ell_2 \in \mathbb{Z}$.
It can be checked that the two contiguous relations in \eqref{dcont} lead to the same recurrence relation
 \[
h(n+1,m,k)=-q^{n+m}h(n,m,k),
 \]
where $n, m, k\in \mathbb{Z}$.

Note that for a vector $x=(x_1, x_2, \ldots, x_r)$,  let $ x\cdot \gamma_i$ denote
the inner product of $x$ and $\gamma_i$, that is,
\[
 x \cdot \gamma_i =x_1\gamma_{i,1}+x_2\gamma_{i,2}+\cdots+x_r\gamma_{i,r}.
\]  For a given multiple theta function $\theta(a)$ of the form \eqref{thetaeq}, we show that $\theta(a)$ satisfies a  contiguous relation of the form \eqref{recmulti}  if and only if there exists a vector $x=(x_1,x_2,\ldots,x_r)\in \mathbb{Q}^r$ such that $\frac{x\cdot \gamma_i}{t_i}$ are integers for $1\leq i \leq m$ and
\begin{equation}\label{rematrix}
v^T=Ax^T,
\end{equation}
where $T$ indicates the transpose of a vector, and  $A$ is an $r\times r$ matrix  defined by
\begin{equation}\label{matrixA0}
A=\left(
  \begin{array}{cccc}
    \frac{\gamma_{1,1}}{t_1} & \frac{\gamma_{2,1}}{t_2} & \cdots & \frac{\gamma_{m,1}}{t_m}  \\[5pt]
    \frac{\gamma_{1,2}}{t_1} &  \frac{\gamma_{2,2}}{t_2}  & \cdots & \frac{\gamma_{m,2}}{t_m} \\[5pt]
    \vdots & \vdots  & \ddots & \vdots \\[5pt]
    \frac{\gamma_{1,r}}{t_1}  &  \frac{\gamma_{2,r}}{t_2}   & \cdots & \frac{\gamma_{m,r}}{t_m} \\
  \end{array}
\right)\left(
  \begin{array}{cccc}
    \gamma_{1,1} & \gamma_{1,2} & \cdots & \gamma_{1,r}  \\[5pt]
    \gamma_{2,1} & \gamma_{2,2} & \cdots & \gamma_{2,r}  \\[5pt]
    \vdots & \vdots &  \ddots & \vdots \\[5pt]
    \gamma_{m,1}  & \gamma_{m,2}  & \cdots & \gamma_{m,r}   \\
  \end{array}
\right).
\end{equation}
By showing that the rank of $A$ is $m$, we can derive  $m$ contiguous relations with linearly independent exponent vectors satisfied by $\theta(a)$.  Denote the exponent vectors of these $m$ contiguous relations by $\nu_1,\nu_2, \ldots, \nu_m$. Moreover, we also verify that any exponent vector of $a$ in the expansion of $\theta(a)$ can be represented as a linear combination of  $\nu_1,\nu_2, \ldots, \nu_m$ with rational coefficients. Then by applying Stanley's Lemma on the fundamental parallelepiped, it leads to that all the coefficients of $\theta(a)$ can be reduced to a finite number of initial values by using the $m$ recurrence relations derived from these $m$ contiguous relations of $\theta(a)$.

In this paper, we shall present an algorithmic approach to the verification of
multiple theta function identities of the form
  \begin{equation}\label{multieq}
\theta_{n+1}(a)=\sum_{k=1}^n c_k \theta_k(a),
\end{equation}
where $a=(a_1,a_2,\ldots, a_r)$, $n\ge 1$, each $c_k$ is a nonzero  complex number,
and for $1\leq k \leq n+1$, each $\theta_k(a)$ is of the form \eqref{thetaeq} which is given by
\begin{equation}\label{thetaak}
\theta_k(a)=a^{\tau_k}\prod_{i=1}^m [(-1)^{\delta_{k,i}}a^{\gamma_{i}^{(k)}}q^{s_{k,i}}; q^{t_{k,i}}]_\infty,
\end{equation}
where $\gamma_{i}^{(k)} \in \mathbb{Z}^r$, $\delta_{k,\,i} =0$ or $1$, $s_{k,\,i} \in  \mathbb{Q}$, $t_{k,\,i}\in \mathbb{Q}^{+}$, $\gamma_1^{(k)}, \gamma_2^{(k)}, \ldots, \gamma_m^{(k)}$ are linearly independent, and $\tau_k$ is a linear combination of $\gamma_1^{(k)}, \gamma_2^{(k)}, \ldots, \gamma_m^{(k)}$ over $\mathbb{Q}$. Furthermore, we assume that $\theta_1(a), \theta_2(a), \ldots, \theta_n(a)$ are linearly independent.
In fact, this restriction occurs often in the literature of theta function identities, and it can be verified by direct computation. For each $1\leq k \leq n+1$,  let $A_k$ be the matrix associated with $\theta_k(a)$ as defined by \eqref{matrixA0}.
 We show that, if
\begin{equation}\label{coniden1}
A_1=A_2=\cdots=A_{n+1},
\end{equation}
there are $m$ contiguous relations with linearly independent exponent vectors satisfied by all of $\theta_1(a), \theta_2(a), \ldots, \theta_{n+1}(a)$, which lead to $m$ recurrence relations on the coefficients of each $\theta_k(a)$. Then by applying Stanley's lemma on the fundamental parallelepiped,  such an identity can be reduced to a finite number of simpler relations on the coefficients of $\theta_1(a), \theta_2(a), \ldots, \theta_{n+1}(a)$.

Let us use the following example to illustrate the steps to verify \eqref{multieq}:
\begin{equation}\label{ideab}
(a,q/a,-b,-q/b\,;q)_\infty+(-a,-q/a,b,q/b\,;q)_\infty
=\frac{2(ab,q^2/ab,aq/b,bq/a;q^2)_\infty}{(q;q^2)_\infty^2}.
\end{equation}
This addition formula can be found in Berndt \cite[p.\,45]{berndt85}.
Clearly, the identity \eqref{ideab} is of the form  (\ref{multieq}), that is,  $r=n=m=2$, and  it can be written as
 \[
 \theta_3=-\theta_1+\frac{2}{(q;q^2)_\infty^2}\theta_2,
 \]
 where
 \begin{align*}
&\theta_1=(-a,-q/a,b,q/b\,;q)_\infty,\\[5pt]
&\theta_2= (ab,q^2/ab,aq/b,bq/a;q^2)_\infty,\\[5pt]
&\theta_3=(a,q/a,-b,-q/b\,;q)_\infty.
\end{align*}
It is easy to see that $\theta_1, \theta_2$ are linearly independent and $A_1=A_2=A_3$, where each $A_i$ is given as \eqref{matrixA0}.  The identity  \eqref{ideab} can be proved via the following steps.

{\noindent \bf Step 1.} For a fixed $1\leq k\leq n+1$, let $\theta_k(a)$ be given as \eqref{thetaak},
where $\gamma_{1}^{(k)},\gamma_{2}^{(k)}, \ldots, \gamma_{m}^{(k)}$ are the exponent vectors in the factors of $\theta_k(a)$.
We first show that if  $\theta_{n+1}(a)$ satisfies a contiguous relation of the form
\begin{equation}\label{conrthetak}
\frac{\theta_{n+1}(aq^x)}{\theta_{n+1}(a)}=\frac{(-1)^{\rho}}{q^{s}a^\nu}
\end{equation}
then $\frac{x\cdot \gamma_i^{(n+1)}}{t_{n+1,i}}$ are integers for $1\leq i\leq m$ and
 \begin{equation}\label{exvcon}
 \nu^T=A_{n+1} x^T,
 \end{equation}
where  $aq^x=(a_1q^{x_1}, a_2q^{x_2}, \ldots, a_rq^{x_r})$, $\rho=0$ or $1$, $s\in \mathbb{Q}$, $v\in \mathbb{Z}^r$ and $x=(x_1,x_2, \ldots, x_r)\in \mathbb{Q}^r$ is a nonzero vector.
By showing that  the rank of the matrix $A_{n+1}$ is equal to $m$, the relation \eqref{exvcon} will lead to $m$ contiguous relations with linearly independent exponent vectors satisfied by $\theta_{n+1}(a)$.

Next, under the assumption that  $A_1=A_2=\cdots=A_{n+1}$, it is always possible to
multiply $x$ by a sufficient large enough integer $N$ such that
 $\frac{Nx\cdot \gamma_i^{(k)}}{t_{k,i}}$ is an integer  for any $1\leq k \leq n+1$ and $1\leq i \leq m$.
Consequently, the relation \eqref{exvcon}  leads  to a contiguous relation with   exponent vector $N\nu$ satisfied by all the the multiple theta functions $\theta_1(a), \theta_2(a), \ldots, \theta_{n+1}(a)$.
This procedure gives $m$ contiguous relations with linearly independent exponent vectors satisfied by all the multiple theta functions $\theta_1(a), \theta_2(a), \ldots, \theta_{n+1}(a)$.
 We shall use $W=\{w_1,w_2,\ldots, w_m\}$  to denote the set of the exponent vectors of these  contiguous relations.

For example, for $\theta_3$ in the identity \eqref{ideab},  we have
\[
A_3= \left(
  \begin{array}{cc}
    1 & 0 \\[5pt]
   0& 1
  \end{array}
\right).
\]
Choosing  $x$ from
\[
B=\{(1,0),(1,1)\},
\]
the relation (\ref{exvcon}) gives two contiguous relations:
\begin{align}
&\frac{\theta_3(aq,b)}{\theta_3(a,b)}=-\frac{1}{a}, \label{conab-1}\\[5pt]
&\frac{\theta_3(aq,bq)}{\theta_3(a,b)}=-\frac{1}{ab}.\label{conab-2}
\end{align}
While  $\theta_1$ and $\theta_2$ also satisfy the contiguous relation \eqref{conab-2},  they do not satisfy the contiguous relation \eqref{conab-1}.  By replacing  $(1,0)$ with  $(2,0)$ in $B$,  we derive the following contiguous relation
satisfied by all of $\theta_1$, $\theta_2$ and $\theta_3$:
 \begin{equation}\label{conab-3}
 \frac{\theta_3(aq^2,b)}{\theta_3(a,b)}=\frac{1}{a^2q}.
 \end{equation}
So we obtain two contiguous relations \eqref{conab-2} and \eqref{conab-3} with linearly independent exponent vectors that are satisfied by  $\theta_1, \theta_2$ and $\theta_3$, where $W=\{w_1,w_2\}=\{(1,1), (2,0)\}$.

{\noindent \bf Step 2.} From the $m$ contiguous relations associated with the vectors in $W$, we can produce $m$ recurrence relations satisfied by the coefficients of each $\theta_k(a)$. Note that for $1\leq k \leq n+1$, $\theta_k(a)$ can be expanded by using Jacobi's triple product identity \eqref{jactri}.
More precisely, $\theta_k(a)$ as defined in \eqref{thetaak} can be expanded as follows:
\begin{align}\label{extheta2}
\theta_k(a)&=a^{\tau_k}\prod_{i=1}^m [(-1)^{\delta_{k,i}}a^{\gamma_{i}^{(k)}}q^{s_{k,i}}; q^{t_{k,i}}]_\infty \nonumber\\[5pt]
&=a^{\tau_k}\prod_{i=1}^m \frac{1}{(q^{t_{k,i}};q^{t_{k,i}})_\infty}\sum_{\ell_i=-\infty}^{\infty} (-1)^{(\delta_{k,i}+1)\ell_i}q^{t_{k,i}{\ell_i\choose 2}}\big(a^{\gamma_{i}^{(k)}}q^{s_{k,i}}\big)^{\ell_i},
\end{align}
which can be written as a multiple sum:
\begin{equation}\label{extheta1}
\theta_k(a)=\sum_{\eta\in \mathbb{Z}^r} h_{k,\eta} a^{\eta},
\end{equation}
where the sum ranges over the vectors $\eta$ such that
\begin{equation}\label{exapeta}
 \eta=\tau_k+\ell_1 \gamma_1^{(k)}+\ell_2\gamma_2^{(k)}+\cdots+\ell_m\gamma_m^{(k)}
\end{equation}
with $\ell_i \in \mathbb{Z}\ (1\leq i \leq m)$. By the assumption that  $\tau_k$ is a linear combination of $\gamma_1^{(k)}, \gamma_2^{(k)}, \ldots, \gamma_m^{(k)}$ over $\mathbb{Q}$, we see that any $\eta$ as given by \eqref{exapeta} can be represented as a linear combination of $\gamma_1^{(k)}, \gamma_2^{(k)}, \ldots,\gamma_m^{(k)}$ with rational coefficients.
 For each $k$,  we also show that  $\eta$ can be expressed as a linear combination of $w_1, w_2, \ldots, w_m$ with rational coefficients.
Using Stanley's Lemma on the fundamental parallelepiped and the recurrence relations derived from these $m$ contiguous relations, we see  that  the coefficients in the expansion of  $\theta_k(a)$ can be reduced to a finite number of initial values with exponent vectors in the set
\[
\Pi_W=\{\lambda_1w_1+\lambda_2w_2+\cdots+\lambda_mw_m \ |\
0\leq\lambda_i <1, \ 1\leq i \leq m\}\cap \mathbb{Z}^r.
\]
 For a given set $W$ consisting of linearly independent integer vectors, we also provide an algorithm to compute $\Pi_W$.

For example, for the contiguous relations \eqref{conab-2} and \eqref{conab-3}, we have
\[
W=\{(1,1), (2,0)\}
\]
and
\begin{align}\label{abpi}
\Pi_W &=\{\lambda_1(1,1)+\lambda_2
(2,0)\ |\ 0\leq\lambda_1,  \lambda_2<1\}\cap
\mathbb{Z}^2=\{(0,0), (1,0)\}.
\end{align}

{{\noindent \bf Step 3.}} Assume that $\Pi_W=\{\beta_1,\beta_2,\ldots, \beta_{d}\}$. For $1\le k \le n+1$,  let $h_{k,\beta_{i}}$ denote the coefficient of $a^{\beta_i}$ in the expansion \eqref{extheta1} of $\theta_k(a)$.   Applying the recurrence relations obtained form the $m$ contiguous relations satisfied by $\theta_k(a)$ and using Stanley's Lemma, we see that  the coefficients $h_{k,\eta}$ of $\theta_k(a)$ can be determined by the coefficients of $a$ with exponent vectors  $\beta_i\in \Pi_W$.  Thus the identity \eqref{multieq} can be proved by verifying the relations
 \begin{equation}\label{relationbeta}
h_{n+1,\beta_i}=\sum_{k=1}^n c_k h_{k,\beta_i}
\end{equation}
for all $\beta_i\in \Pi_W$.
For any exponent vector $\eta$ of $a$ in the expansion of $\theta_k(a)$, it can be uniquely
expressed in the form \eqref{exapeta}, that is,
\[
\eta=\tau_k+\ell_1\gamma_1^{(k)}+\ell_2\gamma_2^{(k)}+\cdots+\ell_m \gamma_m^{(k)},
\]
where $\ell_i \in \mathbb{Z}$.
Then the coefficient $h_{k,\eta}$ of  $\theta_k(a)$ equals
\begin{equation}\label{heta11}
\prod_{i=1}^m \frac{(-1)^{(\delta_{k,i}+1)\ell_i} q^{t_{k,i} {\ell_i \choose 2}+s_{k,i}\ell_i}}{(q^{t_{k,i}};q^{t_{k,i}})_\infty}.
\end{equation}
Thus, to prove (\ref{multieq}), it suffices to verify \eqref{relationbeta} for all the $\beta_i\in \Pi_W$.

For example, for the addition formula \eqref{ideab}, we have $\Pi_W=\{(0,0), (1,0)\}$.
Let  $\beta_1=(0,0)$ and $\beta_2=(1,0)$. By Jacobi's triple product identity \eqref{jactri}, we have
 \begin{align}\label{jacobiexb}
h_{3,\beta_1}&=[a^0b^0]\,\theta_3=[a^0b^0]\, (a,q/a;q)_\infty (-b,-q/b;q)_\infty \nonumber\\[5pt]
&= [a^0b^0]\, \frac{1}{(q;q)_\infty^2} \sum_{n=-\infty}^\infty (-1)^{n} q^{n \choose 2} a^{n} \sum_{n=-\infty}^\infty (-1)^{n} q^{n \choose 2} (-b)^{n}\nonumber \\[5pt]
&=\frac{1}{(q;q)_\infty^2}.
 \end{align}
Similarly, we find that
\[
h_{1,\beta_1}=\frac{1}{(q;q)_\infty^2}, \quad h_{2,\beta_1}=\frac{1}{(q^2;q^2)_\infty^2},
\]
and
\begin{align*}
& h_{1,\beta_2}=\frac{1}{(q;q)_\infty^2},\
& h_{2,\beta_2}=0,\
&& h_{3,\beta_2}=-\frac{1}{(q;q)_\infty^2}.
\end{align*}
Thus the relations
\begin{equation}\label{linearreh1}
h_{3,\beta_i}=-h_{1,\beta_i}+\frac{2}{(q;q^2)_\infty^2}h_{2,\beta_i}
\end{equation}
hold for  $1\leq i \leq 2$. This proves  \eqref{ideab}. \qed

It should be noted that theta functions satisfying the contiguous relation
\begin{equation}\label{onedimcr}
\frac{f(zq^r)}{f(z)}=\frac{1}{z^n q^m}
\end{equation}
with $r,n$ being positive integers and $m$ being a nonnegative integer,
form a vector space of dimension $n$ over $\mathbb{C}$. This is a classical result in  algebraic
geometry (see, for example, \cite[p.\,212, Theorem 1]{shareid96}).
In fact,  for the contiguous relation \eqref{onedimcr}, our approach also gives
\[
\Pi_W=\{\lambda\cdot n \,|\, 0\leq \lambda<1\}=\{0,1,2, \ldots, n-1\}.
\]
This implies that a theta function  satisfying the contiguous relation \eqref{onedimcr} can be determined by the coefficients of $z^0, z^1, \ldots, z^{n-1}$.  Thus an identity on such theta functions  can be reduced to   relations on the coefficients of $z^0, z^1, \ldots, z^{n-1}$. This approach applies
 to many classical theta function identities, such as Watson's quintuple product identity \cite{watson29} and the septuple product identity due to  Farkars and Kra \cite{farkra99}. Our method can be seen as the  systematic generalization of this approach for the verification of identities on theta functions in one variable.

This paper is organized as follows. The objective of Section~2 is to show
that for a   multiple theta function $\theta(a)$ in the form of \eqref{thetaeq},
 there exist $m$ contiguous relations with linearly independent exponent vectors satisfied by $\theta(a)$.
In Section 3, applying Stanley's Lemma on the fundamental parallelepiped, we reduce the coefficients of $\theta(a)$ to a finite number of the initial values.
Section~4 provides a procedure to obtain $m$ contiguous relations with linearly independent exponent vectors satisfied by all the multiple theta functions in an identity of the form \eqref{multieq}. Examples  are given in
Section~5, including the extended Riemann identity and the addition formulae on the Jacobi theta functions.

\section{Contiguous relations} \label{secmul}

Assume that $\theta(a)$ is a multiple theta function in the form of \eqref{thetaeq}, namely,
\[
\theta (a)=a^{\tau} \prod\limits_{i=1}^m [(-1)^{\delta_i}a^{\gamma_{i}}q^{s_i}; q^{t_i}]_\infty,
\]
where $\tau\in \mathbb{Q}^r$, $1\le m\le r$, and for $1\leq i \leq m$,  $\gamma_{i}=(\gamma_{i,1},\gamma_{i,2},\ldots,\gamma_{i,r})\in \mathbb{Z}^r$, $\delta_i =0$ or $1$, $s_i \in  \mathbb{Q}$, $t_i\in \mathbb{Q}^{+}$.
We further assume that $\gamma_1,\gamma_2, \ldots, \gamma_m$ are linearly independent and $\tau$ is a linear combination of $\gamma_1, \gamma_2, \ldots, \gamma_m$ over $\mathbb{Q}$.  For the above
multiple theta function $\theta(a)$, we define the matrix $A(\theta)$, or simply $A$ if no confusion arises, as follows
\begin{equation}\label{matrixA}
A=\left(
  \begin{array}{cccc}
    \frac{\gamma_{1,1}}{t_1} & \frac{\gamma_{2,1}}{t_2} & \cdots & \frac{\gamma_{m,1}}{t_m}  \\[5pt]
    \frac{\gamma_{1,2}}{t_1} &  \frac{\gamma_{2,2}}{t_2}  & \cdots & \frac{\gamma_{m,2}}{t_m} \\[5pt]
    \vdots & \vdots  & \ddots & \vdots \\[5pt]
    \frac{\gamma_{1,r}}{t_1}  &  \frac{\gamma_{2,r}}{t_2}   & \cdots & \frac{\gamma_{m,r}}{t_m} \\
  \end{array}
\right)  \left(
  \begin{array}{cccc}
    \gamma_{1,1} & \gamma_{1,2} & \cdots & \gamma_{1,r}  \\[5pt]
    \gamma_{2,1} & \gamma_{2,2} & \cdots & \gamma_{2,r}  \\[5pt]
    \vdots & \vdots &  \ddots & \vdots \\[5pt]
    \gamma_{m,1}  & \gamma_{m,2}  & \cdots & \gamma_{m,r}   \\
  \end{array}
\right),
\end{equation}
which is an $r\times r$ matrix of rank $m$.

In this section, we provide an algorithm to produce $m$ contiguous relations with linearly independent exponent vectors satisfied by $\theta(a)$ which are of the form \eqref{recmulti}, namely,
\begin{equation}\label{conthetaa}
\frac{\theta(aq^x)}{\theta(a)}=\frac{(-1)^\rho}{q^s a^v},
\end{equation}
where $x=(x_1,x_2, \ldots, x_r)\in \mathbb{Q}^r$, $aq^x=(a_1q^{x_1},a_2q^{x_2}, \ldots, a_rq^{x_r})$, $\rho=0$ or 1, $s\in \mathbb{Q}$ and $v=(v_1,v_2, \ldots, v_r) \in \mathbb{Z}^r$.
Bear in mind that for the purpose of this paper,
we are only concerned with contiguous relations of $\theta(a)$ that are of the form \eqref{conthetaa}.
It can be shown that $\theta(a)$ satisfies a  contiguous relation of the form \eqref{conthetaa} if and only if there exists a vector $x=(x_1,x_2,\ldots,x_r)\in \mathbb{Q}^r$ such that $\frac{x\cdot \gamma_i}{t_i}$ are integers for $1\leq i \leq m$ and
\begin{equation}\label{rematrix}
v^T=Ax^T.
\end{equation}

Under the assumption that $\gamma_1, \gamma_2, \ldots, \gamma_m$ are linearly independent,
 it is easy to see that the rank of $A$ is $m$.
 Using this fact along with the relation \eqref{rematrix}, we can derive $m$  contiguous relations of $\theta(a)$  in the form of (\ref{conthetaa})  with linearly independent exponent vectors. To this end, we need the following lemma.

\begin{lem}\label{lemcon} Let
\[
f(a)=(a^v,q^t/a^v;q^t)_\infty
\]
be a theta function, where $v$ is an integer and $t$ is a positive rational number.
Assume that $f(a)$ has a  contiguous relation of the form
\begin{equation}\label{confa}
f(aq^\alpha)=\frac{(-1)^\rho}{ q^u a^w }f(a),
\end{equation}
where $w\in \mathbb{Z}$, $\rho=0$ or $1$ and $\alpha, u\in \mathbb{Q}$. Then
$v\,|\,w$ and $\alpha v$ is an integer multiple of $t$.
\end{lem}

\pf By means of Jacobi's triple product identity \eqref{jactri}, $f(a)$ can be written as
\begin{equation}\label{expfa}
f(a)=\frac{1}{(q^{t};q^t)_\infty} \sum_{n=-\infty}^\infty (-1)^n q^{t{n\choose 2}}a^{v n}.
\end{equation}
Substituting $a$ with $aq^{\alpha}$ in (\ref{expfa}),
we get an expansion of $f(aq^\alpha)$. Observing that the powers of $a$
in the expansions of both sides of \eqref{confa} are all multiples of $v$, we find  that $v\,|\,w$.
Plugging the expansions of $f(a)$ and $f(aq^{\alpha})$
into the contiguous relation (\ref{confa}) and equating coefficients of
$a^{v n}$, we are led to
 \[
 (-1)^{\rho+n} q^{t {n\choose 2}+\alpha v n}=(-1)^{n+\frac{w}{v}} q^{t {n+\frac{w}{v}\choose 2}-u},
 \]
which gives
 \[
t{n\choose 2}+\alpha v n=t {n+\frac{w}{v}\choose 2}-u,
 \]
and whence
 \[
\Big(t \frac{w}{v}-\alpha v\Big) n+t{\frac{w}{v}\choose 2}-u=0
 \]
 for all $n \in \mathbb{Z}$. It follows that
  $t \frac{w}{v}=\alpha v$.  This completes the proof. \qed

The following theorem gives a necessary and sufficient condition on the existence of
 the contiguous relations of the form \eqref{conthetaa} for a given multiple theta functions of the form \eqref{thetaeq}.

\begin{thm}\label{thmexv}
 Let $\theta(a)$ be a multiple theta function in the form of \eqref{thetaeq} and let
 $A$ be the $r\times r$ matrix given by \eqref{matrixA}. Then if $\theta(a)$ satisfies a contiguous relation of the form
\begin{equation}\label{contigoustheta}
\frac{\theta(aq^x)}{\theta(a)}=\frac{(-1)^\rho}{q^s a^v},
\end{equation}
with $x=(x_1,x_2,\ldots, x_r)\in \mathbb{Q}^r$ being a nonzero vector, $\rho=0$ or $1$, $s \in \mathbb{Q}$ and $v=(v_1,v_2, \ldots, v_r)\in \mathbb{Z}^r$,  then $\frac{x\cdot \gamma_i}{t_i}$ are integers for $1\leq i\leq m$ and
\begin{equation}\label{matrixA1}
v^T=A x^T.
\end{equation}
Conversely, if there exists a nonzero vector $x$ such that $\frac{x\cdot \gamma_i}{t_i}$ are integers for $1\leq i\leq m$, then \eqref{matrixA1} leads to an integer vector $v=(v_1,v_2, \ldots, v_r)$ which corresponds to a contiguous relation of $\theta(a)$ of the form \eqref{contigoustheta}.
\end{thm}

\pf First, assume that  $\theta(a)$ satisfies the contiguous relation \eqref{contigoustheta}, we proceed to show that $\frac{x\cdot \gamma_i}{t_i}$ are integers for $1\leq i\leq m$ and the relation \eqref{matrixA1} holds. Recall that
\[
\theta (a)=a^{\tau} \prod\limits_{i=1}^m [(-1)^{\delta_i}a^{\gamma_{i}}q^{s_i}; q^{t_i}]_\infty,
\]
where  $\gamma_1,\gamma_2, \ldots, \gamma_m$ are linearly independent and $\tau$ is a linear combination of $\gamma_1, \gamma_2, \ldots, \gamma_m$ over $\mathbb{Q}$. Substituting the expression of $\theta(a)$ into \eqref{contigoustheta}, we obtain that
\begin{align*}
\frac{\theta(aq^x)}{\theta(a)}=\frac{a^\tau q^{x\cdot \tau} \prod\limits_{i=1}^m [(-1)^{\delta_i}a^{\gamma_i} q^{x\cdot \gamma_{i} +s_i}; q^{t_i}]_\infty}{a^{\tau} \prod\limits_{i=1}^m [(-1)^{\delta_i}a^{\gamma_{i}}q^{s_i}; q^{t_i}]_\infty}.
\end{align*}
By Lemma \ref{lemcon}, we find that $\frac{x \cdot\gamma_i}{t_i}$ are integers for $1\leq i \leq m$. Let $\ell_i=\frac{x \cdot\gamma_i}{t_i}$ for $1\leq i \leq m$. Then we have
\begin{align}\label{convx}
\frac{\theta(aq^x)}{\theta(a)}
&= q^{x\cdot \tau} \frac{\prod\limits_{i=1}^m [(-1)^{\delta_i}a^{\gamma_i} q^{\ell_i t_i +s_i}; q^{t_i}]_\infty}{ \prod\limits_{i=1}^m [(-1)^{\delta_i}a^{\gamma_{i}}q^{s_i}; q^{t_i}]_\infty}\nonumber \\[5pt]
&=q^{x\cdot\tau} \prod_{i=1}^m \frac{(-1)^{(\delta_i+1) \ell_i}}{q^{\ell_is_i+{ \ell_i\choose 2}t_i}a^{\ell_i\gamma_i}}\nonumber\\[5pt]
&=\frac{(-1)^{\sum_{i=1}^m (\delta_i+1)\ell_i }}{q^{\sum_{i=1}^m\ell_i s_i+{\ell_i\choose 2}t_i- x\cdot \tau} {a}^{\ell_1 \gamma_1+\ell_2
\gamma_2+\cdots +\ell_m\gamma_m}}.
\end{align}
Comparing the above expression with \eqref{contigoustheta}, we get
\begin{align*}
v & =\ell_1 \gamma_1+\ell_2
\gamma_2+\cdots +\ell_m\gamma_m\nonumber\\[5pt]
&=\frac{x\cdot \gamma_1}{t_1} \gamma_1+ \frac{x\cdot \gamma_2}{t_2}\gamma_2+\cdots+\frac{x\cdot \gamma_m}{t_m}\gamma_m,
\end{align*}
which can be rewritten in the form of \eqref{matrixA1}.

Conversely, if there exists a nonzero vector $x=(x_1,x_2,\ldots, x_r)\in \mathbb{Q}^r$ such that $\frac{x\cdot \gamma_i}{t_i}$ are integers for $1\leq i\leq m$, then it follows from \eqref{matrixA1}  that
\begin{align*}
v^T&=Ax^T\\[7pt]
&=\left(
  \begin{array}{cccc}
    \frac{\gamma_{1,1}}{t_1} & \frac{\gamma_{2,1}}{t_2} & \cdots & \frac{\gamma_{m,1}}{t_m}  \\[5pt]
    \frac{\gamma_{1,2}}{t_1} &  \frac{\gamma_{2,2}}{t_2}  & \cdots & \frac{\gamma_{m,2}}{t_m} \\[5pt]
    \vdots & \vdots  & \ddots & \vdots \\[5pt]
    \frac{\gamma_{1,r}}{t_1}  &  \frac{\gamma_{2,r}}{t_2}   & \cdots & \frac{\gamma_{m,r}}{t_m} \\
  \end{array}
\right)  \left(
  \begin{array}{cccc}
    \gamma_{1,1} & \gamma_{1,2} & \cdots & \gamma_{1,r}  \\[5pt]
    \gamma_{2,1} & \gamma_{2,2} & \cdots & \gamma_{2,r}  \\[5pt]
    \vdots & \vdots &  \ddots & \vdots \\[5pt]
    \gamma_{m,1}  & \gamma_{m,2}  & \cdots & \gamma_{m,r}   \\
  \end{array}
\right) x^T\\[5pt]
&=\left(\begin{array}{cccc}
    \frac{\gamma_{1}}{t_1} & \frac{\gamma_{2}}{t_2} & \cdots & \frac{\gamma_{m}}{t_m}
  \end{array} \right) \left(
  \begin{array}{c}
    \gamma_{1}\\
    \gamma_{2} \\
    \vdots  \\
    \gamma_{m}
  \end{array}
\right) x^T\\[5pt]
&=\frac{ x\cdot \gamma_1}{t_1} \gamma_1+\frac{x\cdot \gamma_2 }{t_2} \gamma_2+\cdots+\frac{ x\cdot \gamma_m }{t_m} \gamma_m.
\end{align*}
Under the assumption  that  $\gamma_i \in \mathbb{Z}^r$ and $\frac{x\cdot \gamma_i}{t_i}$ are integers for $1\leq i \leq m$, we see that $v \in \mathbb{Z}^r$. Following the procedure as given in \eqref{convx},
 we deduce that $\theta(a)$ satisfies a contiguous relation of the form \eqref{contigoustheta} with the exponent vector $v$. This completes the proof. \qed

As an example, let us consider the multiple theta function as given by \eqref{examtheta}, namely,
\[
\theta(a,b,c,d)=(ab,q/ab;q)_\infty(a/b,bq/a;q)_\infty (cd,q/cd;q)_\infty (c/d,qd/c;q)_\infty.
\]
So we have  $\gamma_1=(1,1,0,0)$, $\gamma_2=(1,-1,0,0)$, $\gamma_3=(0,0,1,1)$, $\gamma_4=(0,0,1,-1)$ and $t_1=t_2=t_3=t_4=1$, which leads to that
\begin{equation}\label{matrixex}
A=\left(
  \begin{array}{cccc}
    1 &1 & 0 & 0  \\[5pt]
    1 &  -1  & 0 & 0 \\[5pt]
    0 & 0  & 1 & 1 \\[5pt]
    0  &  0  & 1 & -1 \\
  \end{array}
\right) \left(
  \begin{array}{cccc}
    1 & 1 & 0 & 0  \\[5pt]
    1 &  -1  & 0 & 0 \\[5pt]
    0 & 0  & 1 & 1 \\[5pt]
    0  &  0  & 1 & -1 \\
  \end{array}
\right)=\left(
  \begin{array}{cccc}
    2 & 0 & 0 & 0  \\[5pt]
    0 &  2  & 0 & 0 \\[5pt]
    0 & 0  & 2 & 0 \\[5pt]
    0  &  0  & 0 & 2 \\
  \end{array}
\right).
\end{equation}
By the above theorem we see that  for each nonzero rational vector $x=(x_1, x_2, x_3, x_4)$ such that $\frac{x\cdot \gamma_i}{t_i}$ are integers for $1\leq i \leq 4$, there is a contiguous relation of $\theta(a)$ in
the form of
\begin{equation}\label{contex}
\frac{\theta(aq^{x_1},aq^{x_2},aq^{x_3},aq^{x_4})}{\theta(a,b,c,d)}=\frac{(-1)^\rho}{q^s a^{v_1}b^{v_2}c^{v_3}d^{v_4}},
\end{equation}
where $\rho=0$ or $1$, $s \in \mathbb{Q}$, and $v=(v_1,v_2,v_3,v_4)$ such that
\begin{equation}\label{eqex}
v^T=A x^T,
\end{equation}
which is an integer vector.

For example,  $x=(1,0,0,0)$ is a feasible vector subject to
the above conditions.  By \eqref{eqex}, we have $v=(2,0,0,0)$,
which leads to the following contiguous relation:
\begin{equation}\label{exconre}
\frac{\theta(aq,b,c,d)}{\theta(a,b,c,d)}=\frac{1}{a^2}.
\end{equation}

In general, the following theorem shows that it is possible to
construct $m$ contiguous relations with linearly independent exponent vectors satisfied by a multiple theta function of the form  \eqref{thetaeq}.

\begin{thm}\label{corexv}
Let $\theta(a)$ be a multiple theta function of the form \eqref{thetaeq} and $A$ be the matrix given by \eqref{matrixA}. Then there are  $m$ contiguous relations  of the form \eqref{conthetaa} with linearly independent exponent vectors satisfied by $\theta(a)$.
\end{thm}

To prove the above theorem, we need the following lemma.

\begin{lem}\label{lem-rankA}
Let $\theta(a)$ be a multiple theta function of the form \eqref{thetaeq} and let $A$ be the matrix given by \eqref{matrixA}. Then the rank of $A$ is $m$.
\end{lem}

\pf Given the expression \eqref{thetaeq} of $\theta(a)$,  define
 \begin{equation}\label{matrixB}
B=
\left(
  \begin{array}{cccc}
    \frac{\gamma_{1,1}}{t_1} & \frac{\gamma_{2,1}}{t_2} & \cdots & \frac{\gamma_{m,1}}{t_m}  \\[5pt]
    \frac{\gamma_{1,2}}{t_1} &  \frac{\gamma_{2,2}}{t_2}  & \cdots & \frac{\gamma_{m,2}}{t_m} \\[5pt]
    \vdots & \vdots  & \ddots & \vdots \\[5pt]
    \frac{\gamma_{1,r}}{t_1}  &  \frac{\gamma_{2,r}}{t_2}   & \cdots & \frac{\gamma_{m,r}}{t_m} \\
  \end{array}
\right)
\end{equation}
and
\begin{equation}\label{matrixC}
C= \left(
  \begin{array}{cccc}
    \gamma_{1,1} & \gamma_{1,2} & \cdots & \gamma_{1,r}  \\[5pt]
    \gamma_{2,1} & \gamma_{2,2} & \cdots & \gamma_{2,r}  \\[5pt]
    \vdots & \vdots &  \ddots & \vdots \\[5pt]
    \gamma_{m,1}  & \gamma_{m,2}  & \cdots & \gamma_{m,r}   \\
  \end{array}
\right),
\end{equation}
 so that $A=BC$. We follow the common notation
  $r(D)$ to stand for the rank of a matrix $D$.
  Since $1\leq m\leq r$, under the assumption that $\gamma_1, \gamma_2, \ldots, \gamma_m$ are linearly independent,  we see that  $r(B)= r(C)=m$.
Consequently, $r(A) \leq \min\{r(B), r(C)\}=m$. On the other hand, by Sylvester's inequality, we
 deduce that $r(A)\geq r(B)+r(C)-m=m$, and whence $r(A)=m$. This completes the proof. \qed

Now we are ready to prove Theorem \ref{corexv}.

{\noindent {\it Proof of Theorem \ref{corexv}.}} By Theorem \ref{thmexv}, we see that from a vector $x=(x_1,x_2, \ldots, x_r) \in \mathbb{Q}^r$ such that $\frac{x\cdot \gamma_i}{t_i}$ are integers for $1\leq i \leq m$, one can obtain a  contiguous relation of $\theta(a)$ in the form of \eqref{conthetaa} with an integer exponent vector $v$ satisfying
\begin{equation}\label{matrixre}
v^T=A x^T.
 \end{equation}

Next, we show that the equation \eqref{matrixre} gives $m$ contiguous relations of $\theta(a)$
if we choose $x$ out of a set of $r$ linearly independent vectors in $\mathbb{Q}^r$.
To this end, we assume that $\mu_1, \mu_2, \ldots, \mu_r$ are linearly independent vectors in $\mathbb{Q}^r$. Moreover, we may assume without loss of generality that  $\frac{\mu_k \cdot \gamma_i}{t_i}$
is an integer  for each $k$ and $i$.

For each $1\leq k \leq r$, since $\frac{\mu_k \cdot \gamma_i}{t_i}$
are integers for $1\leq i\leq m$, by Theorem \ref{thmexv}, we get
\begin{align}\label{contnuk}
\nu_k^T&=A\mu_k^T
\end{align}
is an integer vector which leads to a contiguous relations of $\theta(a)$ in the form of  \eqref{conthetaa}.

It remains to show that there are $m$ linearly independent vectors among $ \nu_1, \nu_2, \ldots, \nu_r $.
Let $C$ be the $r \times r$ matrix formed by the column vectors $\nu_1, \nu_2, \ldots, \nu_r$, and let
$D$ be the $r \times r$ matrix formed by the column vectors $\mu_1, \mu_2, \ldots, \mu_r$. It follows from \eqref{contnuk} that  $C=AD$.
Since $\mu_1, \mu_2, \ldots, \mu_r$ are linearly independent, that is, $D$ is invertible,
we deduce that $r(C)=r(A)=m$. This means that one can choose
$m$ linearly independent vectors out of $\nu_1, \nu_2, \ldots, \nu_r$.
Hence there are $m$ contiguous relations of $\theta(a)$ with linearly independent exponent vectors. This completes the proof. \qed

For example, for the matrix $A$ as given by \eqref{matrixex}, we have $r(A)=4$. Take four linearly independent vectors $\mu_1=(1,0,0,0)$,  $\mu_2=(0,1,0,0)$, $\mu_3=(0,0,1,0)$, and $\mu_4=(0,0,0,1)$.
Now, $\frac{\mu_k\cdot \gamma_i}{t_i}$ are already integers for all $k$ and $i$.
As shown before, $\mu_1$ leads to the contiguous relation \eqref{exconre}.
Applying equation \eqref{eqex} to $\mu_2, \mu_3$, and $\mu_4$, we are led to
the following three contiguous relations of $\theta(a,b,c,d)$:
\begin{align*}
&\frac{\theta(a,bq,c,d)}{\theta(a,b,c,d)}=\frac{1}{b^2q}, \\[5pt]
&\frac{\theta(a,b,cq,d)}{\theta(a,b,c,d)}=\frac{1}{c^2}, \\[5pt]
&\frac{\theta(a,b,c,dq)}{\theta(a,b,c,d)}=\frac{1}{d^2q}.
\end{align*}
These four contiguous relations have linearly independent exponent vectors: $(2,0,0,0)$, $(0,2,0,0)$, $(0,0,2,0)$ and $(0,0,0,2)$.

Denote by $w_1, w_2, \ldots, w_m$  the $m$ linearly independent vectors among  $\nu_1, \nu_2, \ldots, \nu_r$ as given by Theorem \ref{corexv}, which correspond to $m$ contiguous relations of $\theta(a)$ with linearly independent exponent vectors.  The following theorem shows that any exponent vector in the expansion of $\theta(a)$ can be represented as a linear combination of $w_1, w_2, \ldots, w_m$ with rational coefficients.

\begin{thm}\label{lineareq}
Let $\theta(a)$ be a multiple theta function of the form \eqref{thetaeq} which can be expanded as follows by applying Jacobi's triple product identity,
\begin{equation}\label{extheta}
\theta(a)=\sum_{\eta\in \mathbb{Z}^r} h_{\eta} a^{\eta},
\end{equation}
where
 \begin{equation*}
 \eta=\tau+\ell_1 \gamma_1+\ell_2\gamma_2+\cdots+\ell_m\gamma_m
\end{equation*}
with $\ell_i \in \mathbb{Z}$ for $1\leq i \leq m$.
Then any $\eta$ in the expansion \eqref{extheta} can be represented as a linear combination of $w_1, w_2, \ldots, w_m$ with rational coefficients.
\end{thm}

\pf  Recall that in the definition of $\theta(a)$, $\tau$ is a linear combination of $\gamma_1, \gamma_2, \ldots, \gamma_m$ over $\mathbb{Q}$. Thus any $\eta$ in the expansion \eqref{extheta} can be represented as a linear combination of $\gamma_1, \gamma_2, \ldots, \gamma_m$ with rational coefficients. Suppose that
\begin{equation}
\eta=k_1 \gamma_1+k_2 \gamma_2+\cdots+k_m\gamma_m,
\end{equation}
where $k_i\in \mathbb{Q}$.

To show that $\eta$ can be represented as a linear combination of $w_1, w_2, \ldots, w_m$ with rational coefficients, let us consider the following system of linear equations in the variables $x_1,x_2,\ldots,x_r$ and $\varepsilon$:
\begin{equation}\label{linearsys}
\left\{
  \begin{array}{cc}
    x\cdot\gamma_1-\varepsilon k_1t_1=0, \\[5pt]
    x\cdot\gamma_2-\varepsilon k_2t_2=0, \\[5pt]
 \vdots \\[5pt]
    x\cdot\gamma_m-\varepsilon k_mt_m=0,
  \end{array}
\right.
\end{equation}
where $x=(x_1,x_2,\ldots,x_r)$.
We claim that there is a solution such that $x\in \mathbb{Q}^r$ and $\varepsilon \neq 0$.

 Let $M$ denote the coefficient matrix of \eqref{linearsys}, that is,
 \[ M=
\left(
  \begin{array}{ccccc}
    \gamma_{1,1} & \gamma_{1,2} & \cdots & \gamma_{1,r} & -k_1t_1 \\[5pt]
    \gamma_{2,1} & \gamma_{2,2} & \cdots & \gamma_{2,r} & -k_2t_2 \\[5pt]
    \vdots & \vdots & \vdots & \ddots & \vdots \\[5pt]
    \gamma_{m,1}  & \gamma_{m,2}  & \cdots & \gamma_{m,r}  & -k_mt_m \\
  \end{array}
\right)_{m\times (r+1)}.
\]
 Then \eqref{linearsys} can be expressed as
$M\, (x_1,x_2,\ldots,x_r, \varepsilon)^T=0$.
Since $1\leq m\le r$ and $\gamma_1, \gamma_2, \ldots, \gamma_m$ are linearly independent,  we have $r(M)=m$. Thus the solution space $P$ of \eqref{linearsys} in variables $x_1,x_2,\ldots,x_r,\varepsilon$ has dimension $r+1-m$.

Consider another system of linear equations
\begin{equation}\label{linearsys1}
\left\{
  \begin{array}{cc}
    x\cdot\gamma_1=0, \\[5pt]
    x\cdot\gamma_2=0, \\[5pt]
 \vdots \\[5pt]
    x\cdot\gamma_m=0.
  \end{array}
\right.
\end{equation}
 Let $Q$ denote the vector space of solutions $(x_1,x_2,\ldots,x_r)$ of \eqref{linearsys1}.
 Let $Q'$ be the vector space obtained from $Q$ by
 substituting every vector $(x_1,x_2,\ldots,x_r)\in Q$ with $(x_1,x_2,\ldots,x_r,0)$.  Since $\gamma_1,\gamma_2,\ldots,\gamma_m$ are linearly independent, we see that $\dim Q'=\dim Q=r-m$.
 Clearly, any solution $(x_1, x_2, \ldots,x_r)$ of \eqref{linearsys1} gives rise to a solution $(x_1, x_2, \ldots, x_r, 0)$  of  \eqref{linearsys}, which means that $Q'$ is a subspace of $P$. Since $\dim Q'< \dim P$, there exists a solution $(x_1,x_2,\ldots,x_r,\varepsilon)\in P\setminus Q'$ such that  $\varepsilon \neq 0$.

If $(x_1,x_2, \ldots,x_r,\varepsilon)$ is a solution of \eqref{linearsys}, so is an integer multiple of
$(x_1,x_2,\ldots$, $x_r,$ $\varepsilon)$. Consequently, we may assume without loss of generality that  $\varepsilon k_i$ are integers for all $i$,  and hence
\begin{align}\label{computercon}
\frac{\theta(a_1q^{x_1},a_2q^{x_2},\ldots,a_rq^{x_r})}
{\theta(a_1,a_2,\ldots,a_r)}&= q^{x\cdot \tau} \prod_{i=1}^m
\frac{[(-1)^{\delta_i}a_1^{\gamma_{i,1}}a_2^{\gamma_{i,2}}\cdots
a_r^{\gamma_{i,r}}q^{\varepsilon k_it_i+s_i};
q^{t_i}]_\infty}{[(-1)^{\delta_i}a_1^{\gamma_{i,1}}a_2^{\gamma_{i,2}}\cdots
a_r^{\gamma_{i,r}}q^{s_i}; q^{t_i}]_\infty}\nonumber\\[5pt]
&= q^{x\cdot \tau} \prod_{i=1}^m \frac{(-1)^{(\delta_i+1)\varepsilon k_i}}{q^{\varepsilon s_ik_i+{\varepsilon  k_i\choose 2}t_i}(a_1^{\gamma_{i,1}}a_2^{\gamma_{i,2}}\cdots
a_r^{\gamma_{i,r}})^{\varepsilon k_i}}\nonumber\\[5pt]
&=\frac{(-1)^{\sum_{i=1}^m (\delta_i+1)\varepsilon k_i }}{q^{\sum_{i=1}^m \varepsilon s_ik_i+{\varepsilon  k_i\choose 2}t_i-x\cdot \tau} {a}^{\varepsilon (k_1 \gamma_1+k_2
\gamma_2+\cdots +k_m\gamma_m)}}\nonumber\\[5pt]
&=\frac{(-1)^{\sum_{i=1}^m (\delta_i+1)\varepsilon k_i }}{q^{\sum_{i=1}^m \varepsilon s_ik_i+{\varepsilon  k_i\choose 2}t_i -x\cdot \tau }{a}^{\varepsilon
{\eta}}}.
\end{align}

By Theorem \ref{corexv}, it follows that $(\varepsilon {\eta})^T=A x^T$, where the matrix $A$ is given by \eqref{matrixA}. Following the proof of Theorem \ref{corexv}, by the assumption that $\mu_1, \mu_2, \ldots, \mu_r$ are linearly independent over $\mathbb{Q}^r$, we see that
\[
x=c_1\mu_1+c_2\mu_2+\cdots+c_r\mu_r,
\]
where $c_i\in \mathbb{Q}$ for $1\leq i\leq r$, or equivalently,
 \[
x^T=c_1\mu_1^T+c_2\mu_2^T+\cdots+c_r\mu_r^T.
\]
It follows that
\begin{align*}
(\varepsilon {\eta})^T&=c_1A\mu_1^T+c_2A\mu_2^T+\cdots+c_rA\mu_r^T\\[5pt]
&=c_1\nu_1^T+c_2\nu_2^T+\cdots+c_r\nu_r^T,
\end{align*}
where $\nu_k^T=A\mu_k^T$ as given by \eqref{contnuk}.  Since $w_1, w_2, \ldots, w_m$ are the $m$ linearly independent vectors among $\nu_1, \nu_2, \ldots, \nu_r$,  we deduce that $\varepsilon {\eta}$ can be represented as a linear combination of $w_1, w_2, \ldots, w_m$ with rational coefficients, so can $\eta$ since $\varepsilon \neq 0$.  This completes the proof. \qed

\section{Stanley's Lemma and the initial coefficients}\label{mulss}

In this section, we   apply Stanley's Lemma on the  fundamental parallelepiped to reduce the coefficients of a multiple theta function in the form of \eqref{thetaeq}  to a finite number of initial values.

Recall that a multiple theta function $\theta(a)$ in the form of \eqref{thetaeq} can be expanded as \eqref{extheta},
\begin{equation*}
\theta(a)=\sum_{\eta\in \mathbb{Z}^r} h_{\eta} a^{\eta},
\end{equation*}
where
 \begin{equation*}
 \eta=\tau+\ell_1 \gamma_1+\ell_2\gamma_2+\cdots+\ell_m\gamma_m
\end{equation*}
with $\ell_i \in \mathbb{Z}$ for $1\leq i \leq m$.

By Theorem \ref{corexv}, we can find $m$ contiguous relations with linearly independent exponent vectors satisfied by $\theta(a)$. Denote the exponent vectors of these $m$ contiguous relations by  $w_1,w_2, \dots, w_m$. According to the contiguous relation with exponent vector $w_i$, there is a recurrence relation on the coefficients of $\theta(a)$ in the form
\begin{equation}\label{receta1}
h_\eta=(-1)^{\delta_i}q^{\varphi_i(\eta)}h_{\eta-w_i},
\end{equation}
where $1\leq i \leq m$, $\delta_i=0$ or $1$, and $\varphi_i(\eta)\in \mathbb{Q}$ which is related to $\eta$.
Substituting $\eta$ by $\eta+w_i$, the above recurrence relation can be rewritten as
\begin{equation}\label{receta2}
h_\eta=(-1)^{\delta_i}q^{-\varphi_i(\eta+w_i)}h_{\eta+w_i}.
\end{equation}
 Let $b_i$ be a positive integer for $1\leq i\leq m$. By iterating recurrence relation \eqref{receta1} $b_i$  times, we obtain
 \begin{equation}\label{recbi1}
 h_\eta=(-1)^{b_i\delta_i}q^{\zeta_i(b_i)}h_{\eta-b_iw_i},
 \end{equation}
 where $\zeta_i(b_i)\in \mathbb{Q}$ which is related to $b_i$.
 On the other hand, iterating  \eqref{receta2} $b_i$ times, we are led to
 \[
 h_\eta=(-1)^{b_i\delta_i}q^{\zeta_i(-b_i)}h_{\eta+b_iw_i}.
 \]
 Thus \eqref{recbi1} holds for all integers $b_i$.
For $i$ from $1$ to $m$, by iterating the recurrence relations \eqref{recbi1}, we obtain that the $m$ recurrence relations associated with $w_1, w_2, \ldots, w_m$
can be condensed as follows:
\begin{align}\label{etas}
h_\eta&= (-1)^{\delta} q^{s} h_{\eta-b_1w_1-\cdots-b_mw_m},
\end{align}
where $\delta=0$ or 1, and $s\in \mathbb{Q}$. Using Theorem \ref{lineareq}, we have that $\eta$  can be uniquely expressed as
 \begin{equation}\label{recev}
 {\eta}=g_1w_1+g_2w_2+\cdots+g_mw_m,
 \end{equation}
 where $g_i\in \mathbb{Q}$ for $1\leq i \leq m$.  Applying Stanley's Lemma with $b_i=\lfloor g_i\rfloor$ for $1\leq i \leq m$, we see that $\eta-b_1w_1-\cdots-b_mw_m$ falls into a finite set of initial values. Thus the coefficients of $\theta(a)$ can be determined by a finite set of  initial values.

 For example, let
 \[
 \theta(a)=(a_1a_2,q^2/a_1a_2,a_1q/a_2,a_2q/a_1;q^2)_\infty,
  \]
 where $a=(a_1,a_2)$. The procedure in Theorem \ref{corexv} generates the following two contiguous relations:
\begin{align}
&\frac{\theta(a_1q,a_2q)}{\theta(a_1,a_2)}=-\frac{1}{a_1a_2},\label{crabq}\\[5pt]
&\frac{\theta(a_1q^2,a_2)}{\theta(a_1,a_2)}=\frac{1}{a_1^2q},\label{crabq1}
\end{align}
where the exponent vectors $w_1=(1,1)$ and $w_2=(2,0)$
  are linearly independent.
Let
\[
\theta(a)=\sum_{\eta\in \mathbb{Z}^2}h_{\eta}a^\eta,
\]
where $\eta$ ranges over  linear combinations of the vectors $\gamma_1=(1,1)$ and $\gamma_2=(1,-1)$ with integer coefficients. The contiguous relation \eqref{crabq} leads to the recurrence relation
\begin{align}
h_{\eta}&=-q^{\eta \cdot w_1-2}h_{\eta-w_1}. \label{rec3}
\end{align}
Replacing $\eta$ by $\eta+w_1$, \eqref{rec3} can be rewritten as
\begin{align}
h_{\eta}=-q^{-(\eta+w_1) \cdot w_1+2}h_{\eta+w_1}.\label{rec31}
\end{align}
Let $b_1$ be a positive integer. By iterating recurrence relation \eqref{rec3} $b_1$  times, we obtain
\begin{equation}\label{rec32}
h_{\eta}=(-1)^{b_1} q^{ b_1 \eta\cdot  w_1-{b_1\choose 2}w_1\cdot w_1-2b_1} h_{\eta-b_1w_1}.
\end{equation}
On the other hand, iterating  \eqref{rec31} $b_1$ times, we get
\begin{equation}\label{rec33}
h_{\eta}=(-1)^{b_1} q^{ -b_1 \eta\cdot  w_1-{b_1+1\choose 2}w_1\cdot w_1+2b_1} h_{\eta+b_1w_1}.
\end{equation}
Comparing \eqref{rec32} and \eqref{rec33}, we conclude that \eqref{rec32}  holds for all integers $b_1$.

Similarly, the contiguous relation  \eqref{crabq1} implies that
\begin{align}
h_{\eta}&=q^{b_2\eta \cdot w_2-{b_2+1\choose 2} w_2\cdot w_2+b_2} h_{\eta-b_2w_2},\label{rec4}
\end{align}
where $b_2$ is an integer.
Replacing $\eta$ by $\eta-b_1w_1$ in \eqref{rec4} and combining it with  \eqref{rec32},  we deduce that
\begin{align}\label{heta2}
h_\eta &=(-1)^{b_1} q^{ b_1 \eta\cdot  w_1-{b_1\choose 2}w_1\cdot w_1-2b_1} h_{\eta-b_1w_1}\nonumber\\[5pt]
&=(-1)^{b_1}  q^{\eta\cdot(b_1  w_1+b_2 w_2)-{b_1\choose 2}w_1\cdot w_1-{b_2+1\choose 2}w_2\cdot w_2-b_1b_2 w_1\cdot w_2-2b_1+b_2}h_{\eta-b_1w_1-b_2w_2},
\end{align}
where $b_1,b_2$ are integers.  By Theorem \ref{lineareq},
  we see that $\eta$ can be represented as a linear combination of $w_1$ and $w_2$ with rational coefficients. Let
\begin{equation}\label{etag}
\eta=g_1(1,1)+g_2(2,0),
 \end{equation}
so that $g_1$ and $2 g_2$ are integers. Using the recurrence relation  \eqref{heta2} with $b_1=g_1$ and $b_2=\lfloor g_2\rfloor$, we see that $h_\eta$  can be reduced to the coefficient with the following exponent vector:
 \begin{align}\label{initialeta}
 \eta-b_1w_1-b_2w_2&=g_1(1,1)+g_2(2,0)-b_1(1,1)-b_2(2,0)\nonumber \\[5pt]
 &=(2g_2-2\lfloor g_2\rfloor, 0),
 \end{align}
which equals $(0,0)$ or $(1,0)$ since $2g_2$ is an integer.

The above reduction of  the coefficients of $\theta(a)$ turns out to be the same procedure as given by Stanley's Lemma on the fundamental parallelepiped. A fundamental parallelepiped is defined to be the area generated by linearly independent integer vectors $v_1,v_2,\ldots,v_d$ with coefficients belonging to $[0,1)$, that is,
$$
\Pi=\{\lambda_1v_1+\lambda_2v_2+\cdots+\lambda_dv_d \ |\
0\leq\lambda_i <1 \ \mbox{for} \ 1\leq i \leq d\}.
$$
In the study of nonnegative integer solutions of linear homogeneous diophantine equations, Stanley  \cite[Lemma\,4.5.7 (i)]{stanley97} showed that each element of a simplicial monoid $F$ can be determined by an integer point in the fundamental parallelepiped $\Pi$.
Recall that a simplicial monoid  with quasigenerators $v_1,v_2,\ldots,v_d$ is defined to be the set of integer vectors which can be represented as linear combinations of $v_1,v_2,\ldots,v_d$  with rational coefficients, namely,
\begin{equation}\label{monoidf}
F=\{\gamma \in \mathbb{Z}^r\ |\
\gamma=c_1 v_1+c_2v_2+\cdots+c_dv_d , \  c_i \in \mathbb{Q}\ \mbox {for $1 \leq i \leq d$}\}.
\end{equation}

\begin{lem}[Stanley's Lemma] \label{lemv}
Let $v_1,v_2,\ldots,v_d$ be linearly independent integer vectors  of dimension $r$ over $\mathbb{Q}$ and let $F$ be the simplicial monoid  with quasigenerators $v_1,v_2,\ldots,v_d$.  Then every element  $\gamma\in F$ can be expressed uniquely in the form
\begin{equation}\label{gammav}
\gamma=\beta+b_1v_1+b_2v_2+\cdots+b_dv_d,
\end{equation}
where $\beta\in \Pi\cap \mathbb{Z}^r$ and $b_1,b_2,\ldots,b_d$ are integers. Conversely, any  vector $\gamma$ in the form  of \eqref{gammav} belongs to $F$.
\end{lem}

For example, any integer vector $\eta$ in the form of \eqref{etag}  belongs to the following simplicial monoid:
\[
F=\{\gamma \in \mathbb{Z}^2\ |\
\gamma=c_1 (1,1)+c_2(2,0), \  c_1, c_2 \in \mathbb{Q}\}.
\]
Thus, by   Stanley's Lemma, any vector $\eta \in F$ can be written as
\begin{equation*}
\gamma=\beta+b_1(1,1)+b_2(2,0),
\end{equation*}
where $b_1,b_2 \in \mathbb{Z}$, $\beta\in \Pi\cap \mathbb{Z}^2$, and the fundamental parallelepiped $\Pi$ is given by
\[
\Pi=\{\lambda_1(1,1)+\lambda_2(2,0) \mid 0\leq  \lambda_1, \lambda_2 <1\}.
\]
It is easy to see that $ \Pi\cap \mathbb{Z}^2=\{(0,0), (1,0)\}$,
which is the set of the initial values as given by \eqref{initialeta}.

 By Stanley's Lemma, we have the following algorithm.

\begin{thm}\label{maintheorem}
 Let $\theta(a)$ be a multiple theta function in the form of \eqref{thetaeq}.
 Assume that $\theta(a)$ has $m$ contiguous relations with linearly independent exponent vectors $w_1$, $w_2$,  $\ldots,w_m$. Then any coefficient  $h_{\eta}$  in the expansion \eqref{extheta} of $\theta(a)$ can be determined by the recurrence relations derived from the $m$ contiguous relations combined with  one of the initial values in the set
\begin{equation}\label{htheta}
H_\theta=\{h_{\beta_1},h_{\beta_2},\ldots,
h_{\beta_{d}}\},
\end{equation}
where   $d=|\Pi_W|$, $\beta_1,\beta_2,\ldots,\beta_{d} \in \Pi_W$, and
\begin{equation}\label{pib}
\Pi_W=\{\lambda_{1}w_1+\lambda_2w_2+\cdots+\lambda_{m}w_m\
|\ 0\leq\lambda_i <1,\ 1\leq i \leq m\}\cap \mathbb{Z}^r.
\end{equation}
\end{thm}

\pf  Applying Jacobi's triple product identity \eqref{jactri} to the factors of $\theta(a)$ in \eqref{thetaeq},  we get
 \begin{align}\label{expfi}
f_i(a)&=[(-1)^{\delta_i}a^{\gamma_i}q^{s_i}; q^{t_i}]_\infty\nonumber\\[6pt]
&=\frac{1}{(q^{t_i};q^{t_i})_\infty}
\sum_{k_i=-\infty}^\infty (-1)^{(1+\delta_i)k_i} q^{t_i {k_i \choose 2}+s_ik_i} a^{k_i\gamma_i}.
\end{align}
Hence
  \begin{align}\label{productex}
  \theta(a)&=a^{\tau}\prod_{i=1}^m \frac{1}{(q^{t_i};q^{t_i})_\infty}
\sum_{k_i=-\infty}^\infty (-1)^{(1+\delta_i)k_i} q^{t_i {k_i \choose 2}+s_ik_i} a^{k_i\gamma_i}\nonumber\\[6pt]
&=\sum_{\eta \in \mathbb{Z}^r}\Bigg(\sum_{(k_1, k_2,\ldots,k_m)\in \mathbb{Z}^m} \prod_{i=1}^m \frac{(-1)^{(1+\delta_i)k_i} q^{t_i {k_i \choose 2}+s_ik_i}}{(q^{t_i};q^{t_i})_\infty} \Bigg) a^{\eta},
  \end{align}
where $\eta=(\eta_1,\eta_2,\ldots, \eta_r)$ and the inner sum ranges over $(k_1,k_2,\ldots,k_m)\in \mathbb{Z}^m$ such that
\begin{equation}\label{mnlinear}
\eta=\tau+k_1\gamma_1+k_2\gamma_2+\cdots+k_m\gamma_m.
\end{equation}
Comparing \eqref{extheta} and \eqref{productex}, we obtain that
\begin{equation}\label{heta}
h_\eta =\sum_{(k_1,k_2,\ldots,k_m)} \prod_{i=1}^m \frac{(-1)^{(1+\delta_i)k_i} q^{t_i {k_i \choose 2}+s_ik_i}}{(q^{t_i};q^{t_i})_\infty},
\end{equation}
where the sum ranges over $(k_1,k_2,\ldots,k_m)\in \mathbb{Z}^m$ such that
the relation \eqref{mnlinear} holds.
Since $\gamma_1,\gamma_2,\ldots,\gamma_m$ are linearly independent, we see that for a given exponent vector $\eta$,  there is
 a unique vector $(k_1,k_2,\ldots,k_m) \in \mathbb{Z}^{m}$ that
  satisfies \eqref{mnlinear}. By the assumption of $\tau$, it follows from \eqref{mnlinear}  that $\eta$ is a linear combination of $\gamma_1,\gamma_2,\ldots,\gamma_m$ with rational coefficients.

By Theorem \ref{lineareq}, we see that $\eta$  can be expressed as
 \begin{equation}\label{recev}
 {\eta}=g_1w_1+g_2w_2+\cdots+g_mw_m,
 \end{equation}
 where $g_i\in \mathbb{Q}$.
On the other hand,  Stanley's Lemma yields
\begin{equation}\label{recthe}
\eta=\beta+b_1w_1+b_2w_2+\cdots+b_mw_m,
\end{equation}
where $\beta\in \Pi_W$ and $b_1,b_2,\ldots,b_m$ are integers.

Using the above expressions of $\eta$, we shall establish recurrence
relations on the coefficients of $\theta(a)$.
For $1\leq i \leq m$, assume that $w_i$ corresponds to the following contiguous relation of $\theta(a)$:
\begin{equation}\label{recwi}
\frac{\theta(aq^{\alpha_i})}
{\theta(a)}
=\frac{(-1)^{\delta_i}}{q^{u_i}a^{w_i}},
\end{equation}
where $\alpha_i\in \mathbb{Q}^r$.  Rewriting
\eqref{recwi} as
\begin{equation*}
\theta(a)
={(-1)^{\delta_i}}{q^{u_i}a^{w_i}}{\theta(aq^{\alpha_i})}
\end{equation*}
and equating the coefficients  on both sides, we obtain that for $1\leq i\leq m$,
\begin{equation}\label{receta}
h_\eta=(-1)^{\delta_i}q^{u_i+\alpha_i\cdot(\eta-w_i)}h_{\eta-w_i}.
\end{equation}
Replacing $\eta$ by $\eta+w_i$ in \eqref{receta} gives that for $1\leq i\leq m$,
\begin{equation}\label{receta11}
h_{\eta}=(-1)^{\delta_i}q^{-u_i-\alpha_i\cdot\eta}h_{\eta+w_i}.
\end{equation}
Let $b_i$ be a positive integer for $1\leq i\leq m$.  Iterating  \eqref{receta} $b_i$  times, we obtain
\begin{equation}\label{recbi11}
h_\eta=(-1)^{b_i\delta_i}q^{b_iu_i+b_i\alpha_i\cdot \eta-{b_i+1 \choose 2}\alpha_i \cdot w_i}h_{\eta-b_iw_i}.
\end{equation}
 On the other hand, iterating  \eqref{receta11} $b_i$ times, we arrive at
 \begin{equation}\label{recbi12}
 h_\eta=(-1)^{b_i\delta_i}q^{-b_iu_i-b_i\alpha_i \cdot\eta-{b_i \choose 2}\alpha_i \cdot w_i}h_{\eta+b_iw_i}.
 \end{equation}
Thus  \eqref{recbi11} holds for all integers $b_i$.

 Employing the recurrence relations \eqref{recbi11} repeatedly for $i=1$ to $m$, we find that
\begin{align}\label{ceta}
h_\eta&=(-1)^{\sum_{i=1}^m b_i\delta_i} q^{\sum_{i=1}^m \big(b_iu_i-{b_i+1 \choose 2}\alpha_i\cdot w_i+b_i\big(\alpha_i\cdot\eta-\sum_{j=1}^{i-1}b_j\alpha_i\cdot w_j\big)\big)}
h_{\eta-b_1w_1-\cdots-b_mw_m
}.
\end{align}
Now, applying   \eqref{ceta} with $b_1, b_2, \ldots, b_m$ being determined by Stanley's Lemma as given in \eqref{recthe}, we deduce that
\begin{align} \label{cm}
h_\eta &=(-1)^{\sum_{i=1}^m b_i\delta_i} q^{\sum_{i=1}^m \big(b_iu_i-{b_i+1 \choose 2}\alpha_i\cdot w_i+b_i(\alpha_i\cdot\beta+b_i\alpha_i\cdot w_i+\sum_{j=i+1}^{m}b_j\alpha_i\cdot w_j)\big)}h_{\beta}\nonumber\\[6pt]
&=(-1)^{\sum_{i=1}^m b_i\delta_i} q^{\sum_{i=1}^m \big(\frac{\alpha_i\cdot w_i}{2}b_i^2-\frac{\alpha_i\cdot w_i}{2}b_i + b_iu_i+b_i(\alpha_i\cdot\beta+\sum_{j=i+1}^{m}b_j\alpha_i\cdot w_j)\big)}h_{\beta},
\end{align}
where $\beta$  belongs to $\Pi_W$. Clearly, $\Pi_W$ is  finite.
Thus any coefficient $h_\eta$ in the expansion  of $\theta(a)$ can be  reduced to a coefficient $h_\beta$ with $\beta\in \Pi_W$. This completes the proof.  \qed

For the fundamental parallelepiped $\Pi_W$ of  the linearly independent integer vectors in $W=\{w_1,w_2,\ldots,w_m\}$, as pointed out by Stanley \cite{stanley97},  $\Pi_W$ is a finite set, since it is contained in the intersection of the discrete set $F$ as defined in \eqref{monoidf} with the bounded set of all the  vectors $\lambda_1w_1+\lambda_2w_2+ \cdots+\lambda_mw_m \in \mathbb{R}^m$ with $0 \leq \lambda_i < 1$.  In fact,  $\Pi_W$ can be determined by the following procedure.

Let $w_i=(w_{i,\,1},w_{i,\,2},\ldots, w_{i,\,r})$ for $1\leq i\leq m$.
First, we find two vectors $c=(c_1,c_2,\ldots,c_r)$ and $d=(d_1,d_2,\ldots,d_r)$ such that for any $v=(v_1,v_2, \ldots, v_r)$ in $\Pi_W$, we have
\begin{equation}\label{boundcd}
c_j \leq v_j \leq d_j,
\end{equation}
where $1\leq j \leq r$.
To define $c_j$ and $d_j$ for $1\leq j\leq r$, we assume that
the entries $w_{1,j}, w_{2,j}, \ldots, w_{m,j}$ can be
rearranged as follows:
\[
w'_{1,j}\leq w'_{2,j} \leq \cdots \leq  w'_{s,j}\leq 0 \leq  w'_{s+1,j} \leq \cdots  \leq w'_{m,j},
\]
where  $0 \leq s \leq m$ and set $w'_{0,j}=w'_{m+1,j}=0$.
Let
\begin{align*}
&c_j=\min\{w'_{1,j}+ w'_{2,j}+ \cdots +w'_{s,j}+1, 0\},\nonumber\\[5pt]
&d_j=\max\{w'_{s+1,j} +w'_{s+2,j}+ \cdots + w'_{m,j}-1, 0\},
\end{align*}
and let
\[
V=\{(v_1,v_2,\ldots, v_r)\in \mathbb{Z}^r \mid c_j\leq v_j \leq d_j, 1\leq j \leq r\}.
\]
Obviously, $\Pi_W \subseteq V$.  For any $v\in V$, consider the following system of linear equations
in $x_1,x_2, \ldots, x_m$:
\begin{equation}\label{liearwb}
x_1w_1+x_2w_2+ \cdots+x_mw_m=v.
\end{equation}
Since $w_1, w_2, \ldots, w_m$ are linearly independent,  \eqref{liearwb} either has no solution or has a unique solution. If $0\leq x_i < 1$ for $1\leq i \leq m$, then $v\in  \Pi_W$;
otherwise, $v\not\in \Pi_W$.

Note that there are various ways to check whether a system of linear equations \eqref{liearwb} has a solution  $(x_1,x_2, \ldots, x_m)$ with $0\leq x_i <1$ for each $i$, such as Collins' Cylindrical Algebraic Decomposition (CAD) algorithm \cite{co75,ka11} and integer linear programming \cite[Sections 12.2 and 13.4]{sc99}.

\section{Contiguous relations for multiple theta function identities}

In this section, we are concerned with identities on multiple theta function identities of the form \eqref{multieq}, namely,
  \begin{equation}\label{multieq1}
\theta_{n+1}(a)=\sum_{k=1}^n c_k \theta_k(a),
\end{equation}
where each $c_k$ is a nonzero constant, and for  $1\leq k \leq n+1$, $\theta_k(a)$ is given by \eqref{thetaak}.
Furthermore, we assume that $\theta_1(a), \theta_2(a), \ldots, \theta_{n}(a)$  are linearly independent.
Indeed, it is not difficult to check the linear independence by direct computation.
However, as will be seen, for the purpose of verifying an identity, this step may be skipped.

For  $1\leq k \leq n+1$, denote by $A_k$ the $r\times r$ matrix associated with $\theta_k(a)$ as given in \eqref{matrixA}, namely,
\begin{equation}\label{matrixak}
A_k=\left(
  \begin{array}{cccc}
    \frac{\gamma_{1,1}^{(k)}}{t_{k,1}} & \frac{\gamma_{2,1}^{(k)}}{t_{k,2}} & \cdots & \frac{\gamma_{m,1}^{(k)}}{t_{k,m}}  \\[5pt]
    \frac{\gamma_{1,2}^{(k)}}{t_{k,1}} &  \frac{\gamma_{2,2}^{(k)}}{t_{k,2}}  & \cdots & \frac{\gamma_{m,2}^{(k)}}{t_{k,m}} \\[5pt]
    \vdots & \vdots  & \ddots & \vdots \\[5pt]
    \frac{\gamma_{1,r}^{(k)}}{t_{k,1}}  &  \frac{\gamma_{2,r}^{(k)}}{t_{k,2}}   & \cdots & \frac{\gamma_{m,r}^{(k)}}{t_{k,m}} \\
  \end{array}
\right) \left(
  \begin{array}{cccc}
    \gamma_{1,1}^{(k)} & \gamma_{1,2}^{(k)} & \cdots & \gamma_{1,r}^{(k)}  \\[5pt]
    \gamma_{2,1}^{(k)} & \gamma_{2,2}^{(k)} & \cdots & \gamma_{2,r}^{(k)}  \\[5pt]
    \vdots & \vdots &  \ddots & \vdots \\[5pt]
    \gamma_{m,1}^{(k)}  & \gamma_{m,2}^{(k)}  & \cdots & \gamma_{m,r}^{(k)}   \\
  \end{array}
\right).
\end{equation}
For the identity \eqref{multieq1}, it is often the case that
\begin{equation}\label{coniden}
A_1=A_2=\cdots=A_{n+1},
\end{equation}
which ensures that there exist $m$ contiguous relations with linearly independent exponent vectors satisfied by all the multiple theta functions in the identity \eqref{multieq1}.

\begin{thm}\label{vequals}
If a  multiple theta function identity of the form \eqref{multieq1} satisfies that
$A_1=A_2=\cdots= A_{n+1}$, then there exist $m$ contiguous relations with linearly independent exponent vectors satisfied by all the multiple theta functions $\theta_1(a)$, $\theta_2(a)$, $\ldots$, $\theta_{n+1}(a)$.
\end{thm}
\pf First by Theorem \ref{corexv}, we see that there are $m$ contiguous relations with linearly independent exponent vectors satisfied by $\theta_{n+1}(a)$. We proceed to show that from a contiguous relation of $\theta_{n+1}(a)$ with exponent vector $\nu$, one can construct a contiguous relation satisfied by all of $\theta_1(a), \theta_2(a), \ldots, \theta_{n+1}(a)$ with exponent vector $N\nu$ for some positive integer $N$.
This leads to $m$  contiguous relations with linearly independent exponent vectors satisfied by all the multiple theta functions $\theta_1(a), \theta_2(a), \ldots, \theta_{n+1}(a)$.

 Assume that we have the following contiguous relation for $\theta_{n+1}(a)$:
\begin{equation}\label{multiplen}
\frac{\theta_{n+1}(aq^{x})}{\theta_{n+1}(a)}=\frac{(-1)^{\delta}}{q^{s} a^{\nu}},
\end{equation}
where $x=(x_1,x_2,\ldots,x_r)\in \mathbb{Q}^r$ is nonzero, $\nu\in \mathbb{Z}^r$, $\delta=0$ or 1 and $s \in \mathbb{Q}$.  By Theorem \ref{thmexv}, this  contiguous relation of $\theta_{n+1}(a)$
 implies that $\frac{x\cdot \gamma_{i}^{(n+1)}}{t_{n+1,i}}$ are integers for $1\leq i \leq m$, and
\begin{equation}\label{nux}
\nu^T=A_{n+1} x^T.
\end{equation}

We now proceed to construct a contiguous relation with the exponent vector being
a multiple of $\nu$  that is satisfied by all the multiple theta functions $\theta_1(a), \theta_2(a),\ldots, \theta_{n+1}(a)$. To this end, we first determine a positive integer $N$ such that
\begin{equation}\label{integer}
\frac{Nx\cdot\gamma_{1}^{(k)}}{t_{k,1}}, \frac{Nx\cdot\gamma_{2}^{(k)}}{t_{k,2}}, \ldots, \frac{Nx\cdot\gamma_{m}^{(k)}}{t_{k,m}}
\end{equation}
are integers for $1\leq k \leq n+1$.  Clearly, $\theta_{n+1}(a)$ satisfies the following contiguous relation:
 \begin{equation}\label{multiplec1}
\frac{\theta(aq^{Nx})}{\theta(a)}=\frac{(-1)^{\rho}}{q^{u} a^{N\nu}},
\end{equation}
where $\rho=0$ or $1$, and $u\in \mathbb{Q}$.
Then we shall show that all the multiple theta functions $\theta_{1}(a), \theta_2(a), \ldots, \theta_n(a)$ also satisfy the contiguous relation \eqref{multiplec1}.

Consider a multiple theta function $\theta_k(a)$, where  $1\leq k \leq n$.
We claim that $\theta_k(a)$ satisfies a contiguous relation with exponent vector $N\nu$.
We further prove that this contiguous relation is indeed the same as \eqref{multiplec1}.
From the relation \eqref{nux}, it is direct to see that
\[
(N\nu)^T=A_{n+1}\,(Nx)^T.
\]
Then by the assumption that $A_1=A_2=\cdots=A_{n+1}$, we derive that
\begin{equation}\label{nuak}
(N\nu)^T=A_{k}\,(Nx)^T.
\end{equation}
Since $\frac{Nx \cdot \gamma_i^{(k)}}{t_{k,i}}$ is an integer for any $1\leq i\leq m$,
  by Theorem \ref{thmexv},  the above condition ensures that there exists a contiguous relation of
 $\theta_k(a)$ of the following form:
\begin{equation}\label{multipleck}
\frac{\theta_k(aq^{Nx})}{\theta_k(a)}=\frac{(-1)^{\rho_k}}{q^{u_k} a^{N\nu}},
\end{equation}
where $\rho_k=0$ or $1$, and $u_k\in \mathbb{Q}$.

Next we show that if the identity (\ref{multieq1}) holds, then the numbers $(-1)^{\rho_k}$ and $q^{u_k}$ can be uniquely determined.
As a result, the contiguous relation (\ref{multipleck}) satisfied by $\theta_k(a)$
turns out to be exactly the same as the contiguous (\ref{multiplec1}) satisfied by $\theta_{n+1}(a)$.
In doing so, let $d=\frac{(-1)^{\rho}}{q^{u}}$ and $d_k=\frac{(-1)^{\rho_k}}{q^{u_k}}$. Then the contiguous relation  (\ref{multiplec1}) and (\ref{multipleck}) can be written as follows
\begin{equation}\label{multiplec2}
\frac{\theta_{n+1}(aq^{Nx})}{\theta_{n+1}(a)}=\frac{d}{a^{N\nu}},
\end{equation}
and for $1\leq k \leq n$,
\begin{equation}\label{multipleck2}
\frac{\theta_k(aq^{Nx})}{\theta_k(a)}=\frac{d_k}{a^{N\nu}}.
\end{equation}

Substituting $a$ with $aq^{Nx}$ in  \eqref{multieq1}, we get
\begin{align*}
 \theta_{n+1}(aq^{Nx})=\sum_{k=1}^{n} c_k \theta_{k}(aq^{Nx}).
\end{align*}
Applying   (\ref{multiplec2}) and (\ref{multipleck2})
for $1\leq k \leq n$, we find that
\begin{align*}
\frac{d}{ a^{N \nu}} \theta_{n+1}(a)=\sum_{k=1}^{n} c_k \frac{d_k}{ a^{N\nu}} \theta_{k}(a),
\end{align*}
which simplifies to
\begin{align*}
d\, \theta_{n+1}(a)=\sum_{k=1}^{n} c_k d_k \theta_{k}(a).
\end{align*}
Invoking the identity \eqref{multieq1}, we deduce that
\begin{align*}
d\sum_{k=1}^{n} c_k \theta_k(a)=\sum_{k=1}^{n} c_k d_k \theta_{k}(a),
\end{align*}
so that
\begin{align*}
 \sum_{k=1}^{n}c_k (d-d_k) \theta_k(a)=0.
\end{align*}
Since $\theta_1(a), \theta_2(a), \ldots, \theta_{n}(a)$ are linearly independent and
$c_1, c_2, \ldots, c_n$ are all nonzero,
 we conclude that
\[
d=d_1=d_2=\cdots=d_{n}.
\]
Therefore $\rho_k=\rho$, and $u_k=u$ for $1\leq k \leq n$.
This implies that all of $\theta_1(a), \theta_2(a), \ldots$ , $\theta_{n}(a)$ satisfy the contiguous relation \eqref{multiplec1}, and hence the proof is complete. \qed

We remark that when we attempt to prove an identity in the form of (\ref{multieq1}),
we do not have to verify the linear independence of $\theta_1(a), \theta_2(a), \ldots, \theta_n(a)$ as the first step. Instead, we may try to derive the contiguous relations (\ref{multipleck2}) for $1\leq k \leq n$,
and one should expect that $d=d_1=d_2=\cdots=d_n$. If it is indeed the case, then we may
skip the step of justifying the linear independence of $\theta_1(a), \theta_2(a), \ldots, \theta_n(a)$.

\section{Examples}\label{secmulide}

In this section, we give examples to demonstrate how to prove multiple theta
function identities by using our approach.

 As the first example, let us consider Riemann's addition formula (see Krattenthaler \cite{Kra05}).
Weierstrass showed that it is equivalent to an identity on sigma functions (see Whittaker and Watson \cite[p.\,451]{whiwat27}). Note that the addition formulas for theta functions play an important role in the theory of elliptic functions (see, for example, \cite{berndt85,shen94, whiwat27}).

\begin{exam} We have
\begin{align}\label{bailey1}
&\theta(xy)\theta(x/y)\theta(uv)\theta(u/v)-\theta(xv)\theta(x/v)\theta(uy)\theta(u/y)\nonumber\\[5pt]
&\qquad=\frac{u}{y}\theta(yv)\theta(y/v)\theta(xu)\theta(x/u).
\end{align}
\end{exam}

\pf Denote the  multiple theta functions in (\ref{bailey1})  by $\theta_1,\, \theta_2$ and $\theta_3$, namely,
\begin{align*}
\theta_1&=(xy,q/xy,x/y,qy/x,uv,q/uv, u/v, qv/u;q)_\infty,\\[5pt]
\theta_2&=(xv,q/xv,x/v,qv/x,uy,q/uy,u/y,qy/u;q)_\infty,\\[5pt]
\theta_3&= \frac{u}{y}(yv,q/yv,y/v,qv/y,xu,q/xu,x/u,qu/x;q)_\infty.
\end{align*}
Then  \eqref{bailey1} can be written as
\begin{equation}\label{bailey2}
\theta_1=\theta_2+\theta_3.
\end{equation}

By the definition \eqref{matrixak}, for each $\theta_k$,
\[
A_k=\left(
  \begin{array}{cccc}
    2 & 0 & 0 & 0  \\[5pt]
    0 &  2  & 0 & 0 \\[5pt]
    0 & 0  & 2 & 0 \\[5pt]
    0  &  0  & 0 & 2 \\
  \end{array}
\right).
\]
Choosing the vector $x$ from
\[
B=\{(1,0,0,0),(0,1,0,0), (0,0,1,0), (0,0,0,1)\}
\]
and   applying Theorem \ref{thmexv} and  Theorem  \ref{corexv}, we find that $\theta_1$ satisfies the following four  contiguous relations:
\begin{align}\label{conabdef}
&\frac{\theta_1(xq,y,u,v)}{\theta_1(x,y,u,v)}=\frac{1}{x^2}, \nonumber\\[5pt]
&\frac{\theta_1(x,yq,u,v)}{\theta_1(x,y,u,v)}=\frac{1}{y^2q},\nonumber\\[5pt]
&\frac{\theta_1(x,y,uq,v)}{\theta_1(x,y,u,v)}=\frac{1}{u^2}, \nonumber\\[5pt]
&\frac{\theta_1(x,y,u,vq)}{\theta_1(x,y,u,v)}=\frac{1}{v^2q}.
\end{align}
It easily checked  that
 $\theta_2$ and $\theta_3$ also satisfy the above contiguous relations.
Now, we have
\[W=\{(2,0,0,0),
(0,2,0,0),(0,0,2,0),(0,0,0,2)\}
\]
and
\begin{align*}
\Pi_W &=\{\lambda_1(2,0,0,0)+\lambda_2
(0,2,0,0)+\lambda_3(0,0,2,0)\\[5pt]
&\qquad\  +\lambda_4(0,0,0,2) \ |\ 0\leq\lambda_i <1, 1\le i\le 4\}\cap
\mathbb{Z}^4\\[5pt]
&=\{(a_1,a_2,a_3,a_4)\ |\ a_i=0\ \mbox{or}\ 1, 1\le i\le 4\}.
\end{align*}
Because of the symmetries in the parameters $x,y,u$ and $v$, we have only to show that identity \eqref{bailey1} holds for the terms with exponent vectors in
\begin{equation}\label{setww}
\Pi_W'=\{(0,0,0,0), (1,0,0,0), (1,1,0,0), (1,1,1,0), (1,1,1,1)\}.
\end{equation}

Denote the vectors in $\Pi_W'$  by $\beta_1=(0,0,0,0)$, $\beta_2=(1,0,0,0)$, $\beta_3=(1,1,0,0)$, $\beta_4=(1,1,1,0)$ and $\beta_5=(1,1,1,1)$. Then to prove the identity \eqref{bailey1}, it is sufficient to verify the relations
\begin{equation}\label{linearreh}
h_{3,\beta_i}=h_{1,\beta_i}-h_{2,\beta_i},
\end{equation}
for $1\leq i \leq 5$.
Consider the case $k=1$ and $i=1$. In this case, $h_{1,\beta_1}=[x^0y^0u^0v^0]\,\theta_1(a)$. By Jacobi's triple product identity \eqref{jactri}, we have
 \begin{align}\label{jacobiexb}
 [x^0y^0u^0v^0]\,\theta_1(a)&=[x^0y^0u^0v^0]\, (xy,q/xy;q)_\infty (x/y,qy/x;q)_\infty (uv,q/uv;q)_\infty (u/v,qv/u;q)_\infty \nonumber\\[5pt]
&= [x^0y^0u^0v^0]\, \frac{1}{(q;q)_\infty^4} \sum_{n=-\infty}^\infty (-1)^{n} q^{n \choose 2} (xy)^{n} \sum_{n=-\infty}^\infty (-1)^{n} q^{n \choose 2} (x/y)^{n}\nonumber\\[5pt]
&\qquad\qquad\qquad\qquad\quad\times \sum_{n=-\infty}^\infty (-1)^{n} q^{n \choose 2} (uv)^{n} \sum_{n=-\infty}^\infty (-1)^{n} q^{n \choose 2} (u/v)^{n}\nonumber\\[5pt]
&=\frac{1}{(q;q)_\infty^4}\ \sum_{(n_1,n_2,n_3,n_4)\in \mathbb{Z}^4} (-1)^{n_1+n_2+n_3+n_4} q^{{n_1\choose 2}+{n_2\choose 2}+{n_3\choose 2}+{n_4\choose 2}},
 \end{align}
 where the summation in \eqref{jacobiexb}  ranges over integer vectors $(n_1,n_2,n_3,n_4)$ such that
\[
\left\{
  \begin{array}{ll}
    n_1+n_2=0,  \\[4pt]
    n_1-n_2=0,  \\[4pt]
    n_3+n_4=0,  \\[4pt]
    n_3-n_4=0.
  \end{array}
\right.
\]
Solving the above equations, we get
\[
h_{1,\beta_1}=\frac{1}{(q;q)_\infty^4}.
\]
Similarly, we find that
\[
h_{2,\beta_1}=\frac{1}{(q;q)_\infty^4}, \quad h_{3,\beta_1}=0,\\[5pt]
\]
and hence relation \eqref{linearreh} holds for $i=1$.
 For $i=2,3,4, 5$, we obtain
\begin{align*}
&h_{1,\beta_2}=0, && h_{2,\beta_2}=0, && h_{3,\beta_2}=0,\\[5pt]
&h_{1,\beta_3}=\frac{-1}{(q;q)_\infty^4}, && h_{2,\beta_3}=0, &&  h_{3,\beta_3}=\frac{-1}{(q;q)_\infty^4},\\[5pt]
&h_{1,\beta_4}=0, && h_{2,\beta_4}=0,  &&
h_{3,\beta_4}=0,\\[5pt]
&h_{1,\beta_5}=\frac{1}{(q;q)_\infty^4}, && h_{2,\beta_5}=\frac{1}{(q;q)_\infty^4}, && h_{3,\beta_5}=0.
\end{align*}
Thus (\ref{linearreh}) is true for $1\leq i\leq 5$. This proves   (\ref{bailey1}). \qed

The next example is concerned with the extended Riemann identity on theta functions due to Malekar and Bhate   \cite[Theorem 3.1]{MalBha12}.

 \begin{exam} We have
\begin{align}\label{exriemann}
&{\small{ 4qxyuv[-q^2x^2,-q^2y^2,-q^2u^2,-q^2v^2;q^2]_\infty +4q^{\frac{1}{2}}xy[-q^2x^2,-q^2y^2,-qu^2,-qv^2;q^2]_\infty}} \nonumber\\[5pt]
&\qquad+4q^{\frac{1}{2}}uv[-qx^2,-qy^2,-q^2u^2,-q^2v^2;q^2]_\infty +4[-qx^2,-qy^2,-qu^2,-qv^2;q^2]_\infty\nonumber\\[5pt]
&\quad  = \frac{(q^{\frac{1}{2}};q^\frac{1}{2})_\infty^4}{(q^2;q^2)_\infty^4}\Big([-q^{\frac{1}{4}}x,-q^{\frac{1}{4}}y,
-q^{\frac{1}{4}}u,-q^{\frac{1}{4}}v;q^{\frac{1}{2}}]_\infty+[q^{\frac{1}{4}}x,q^{\frac{1}{4}}y,
-q^{\frac{1}{4}}u,-q^{\frac{1}{4}}v;q^{\frac{1}{2}}]_\infty\nonumber\\[5pt]
&\quad\quad\qquad+[-q^{\frac{1}{4}}x,-q^{\frac{1}{4}}y,
q^{\frac{1}{4}}u,q^{\frac{1}{4}}v;q^{\frac{1}{2}}]_\infty
+[q^{\frac{1}{4}}x,q^{\frac{1}{4}}y,
q^{\frac{1}{4}}u,q^{\frac{1}{4}}v;q^{\frac{1}{2}}]_\infty\Big).
\end{align}
\end{exam}

\pf  Denote the multiple theta functions in  \eqref{exriemann} by $\theta_k(x,y,u,v)$ for $1\leq k \leq 8$, namely,
 \begin{align*}
 &\theta_1=xyuv[-q^2x^2,-q^2y^2,-q^2u^2,-q^2v^2;q^2]_\infty, && \theta_2=xy[-q^2x^2,-q^2y^2,-qu^2,-qv^2;q^2]_\infty,\\[5pt]
 &\theta_3=uv[-qx^2,-qy^2,-q^2u^2,-q^2v^2;q^2]_\infty, && \theta_4=[-qx^2,-qy^2,-qu^2,-qv^2;q^2]_\infty,\\[5pt]
 &\theta_5=[-q^{\frac{1}{4}}x,-q^{\frac{1}{4}}y,
-q^{\frac{1}{4}}u,-q^{\frac{1}{4}}v;q^{\frac{1}{2}}]_\infty,  && \theta_6=[q^{\frac{1}{4}}x,q^{\frac{1}{4}}y,
-q^{\frac{1}{4}}u,-q^{\frac{1}{4}}v;q^{\frac{1}{2}}]_\infty,\\[5pt]
&\theta_7=[-q^{\frac{1}{4}}x,-q^{\frac{1}{4}}y,
q^{\frac{1}{4}}u,q^{\frac{1}{4}}v;q^{\frac{1}{2}}]_\infty, &&
\theta_8=[q^{\frac{1}{4}}x,q^{\frac{1}{4}}y,
q^{\frac{1}{4}}u,q^{\frac{1}{4}}v;q^{\frac{1}{2}}]_\infty.
 \end{align*}
 Hence  \eqref{exriemann}  takes the form:
\begin{equation}\label{exriemann8}
\theta_8(x,y,u,v)=\sum_{k=1}^7 c_k\theta_k(x,y,u,v),
\end{equation}
where
\[
c_1= 4q \frac{(q^2;q^2)_\infty^4}{(q^{\frac{1}{2}};q^\frac{1}{2})_\infty^4}, \quad c_2=c_3=4q^{\frac{1}{2}}\frac{(q^2;q^2)_\infty^4}{(q^{\frac{1}{2}};q^\frac{1}{2})_\infty^4},\quad c_4=4\frac{(q^2;q^2)_\infty^4}{(q^{\frac{1}{2}};q^\frac{1}{2})_\infty^4},\quad c_5=c_6=c_7=-1.
\]
 We find that
\begin{align*}
&\frac{\theta_k(xq,y,u,v)}{\theta_k(x,y,u,v)}=\frac{1}{x^2q}, &
&\frac{\theta_k(x,yq,u,v)}{\theta_k(x,y,u,v)}=\frac{1}{y^2q},\\[6pt]
&\frac{\theta_k(x,y,uq,v)}{\theta_k(x,y,u,v)}=\frac{1}{u^2q}, &
&\frac{\theta_k(x,y,u,vq)}{\theta_k(x,y,u,v)}=\frac{1}{v^2q},
\end{align*}
for $k=1, 2, \ldots, 8$. The set of  exponent vectors of the above contiguous relations is
\[
W=\{(2,0,0,0),(0,2,0,0),(0,0,2,0),(0,0,0,2)\}.
 \]
The fundamental parallelepiped of $W$ equals
\begin{align*}
\Pi_W &=\{\lambda_1(2,0,0,0)+\lambda_2
(0,2,0,0)+\lambda_3(0,0,2,0)\\[5pt]
&\qquad +\lambda_4(0,0,0,2) \ |\ 0\leq\lambda_i <1, 1\le i\le 4\}\cap
\mathbb{Z}^4\\[5pt]
&=\{(a_1,a_2,a_3,a_4)\ |\ a_i=0\ \mbox{or}\ 1, 1\le i\le 4\}.
\end{align*}
In view of the symmetry of $x, y, u$ and $v$,  the coefficients of $\theta_k(x,y,u,v)$ are determined by the coefficients with  exponent vectors in
\[
\Pi_W'=\{(0,0,0,0), (1,0,0,0), (1,1,0,0), (1,0,1,0), (1,1,1,0), (1,1,1,1)\}.
\]
Thus  \eqref{exriemann8} can be reduced to equalities on
the coefficients with exponent vectors in $\Pi_W'$. In fact,
for each vector in $\Pi_W'$, the required equality
 immediately follows from Jacobi's triple product identity.
For instance,  consider the coefficients with exponent vector  $(0,0,0,0)$.
By Jacobi's triple product identity, we see that
 \[
 [x^0y^0u^0v^0]\,\theta_8(x,y,u,v)=\frac{1}{(q^{\frac{1}{2}};q^{\frac{1}{2}})_\infty^4}
 \]
and
\begin{align*}
[x^0 & y^0u^0v^0]\, \sum_{k=1}^7 c_k \theta_k(x,y,u,v) =[x^0y^0u^0v^0]\, \sum_{k=4}^7 c_k \theta_k(x,y,u,v)\\[5pt]
&=\frac{4}{(q^{\frac{1}{2}};q^{\frac{1}{2}})_\infty^4}-\frac{1}{(q^{\frac{1}{2}};q^{\frac{1}{2}})_\infty^4}-\frac{1}{(q^{\frac{1}{2}};q^{\frac{1}{2}})_\infty^4}
-\frac{1}{(q^{\frac{1}{2}};q^{\frac{1}{2}})_\infty^4}\\[5pt]
&=\frac{1}{(q^{\frac{1}{2}};q^{\frac{1}{2}})_\infty^4}.
\end{align*}
Similarly, we have
\[
[x^1y^0u^1v^0]\,\theta_8(x,y,u,v)=[x^1 y^0u^1v^0]\, \sum_{k=1}^7 c_k \theta_k(x,y,u,v)=\frac{q^{\frac{1}{2}}}{(q^{\frac{1}{2}};q^{\frac{1}{2}})_\infty^4}.
\]
This completes the proof.   \qed

We now look at an addition formula  on Jacobi theta functions due to Whittaker and Watson \cite[Chapter XXI]{whiwat27}.
Recall that the four classical Jacobi theta functions are defined by
\begin{align*}
\vartheta_1(z)&=\sum_{n=-\infty}^\infty (-1)^{n-\frac{1}{2}} q^{(n+\frac{1}{2})^2}
e^{(2n+1)iz}, \\[5pt]
\vartheta_2(z)&=\sum_{n=-\infty}^\infty  q^{(n+\frac{1}{2})^2}
e^{(2n+1)iz},\\[5pt]
\vartheta_3(z)&=\sum_{n=-\infty}^\infty  q^{n^2}
e^{2niz},\\[5pt]
\vartheta_4(z)&=\sum_{n=-\infty}^\infty (-1)^{n} q^{n^2} e^{2niz}.
\end{align*}

\begin{exam} \label{ex63} We have
\begin{align}\label{add5}
&\vartheta_1(w)\vartheta_1(x)\vartheta_1(y)\vartheta_1(z)
+\vartheta_2(w)\vartheta_2(x)\vartheta_2(y)\vartheta_2(z)\nonumber\\[4pt]
&\quad=\vartheta_1(w')\vartheta_1(x')\vartheta_1(y')\vartheta_1(z')
+\vartheta_2(w')\vartheta_2(x')\vartheta_2(y')
\vartheta_2(z'),
\end{align}
where
\begin{align}\label{wxyz}
&2w'=-w+x+y+z,\qquad
2x'=w-x+y+z,\nonumber\\[3pt]
&2y'=w+x-y+z,\qquad\quad
2z'=w+x+y-z.
\end{align}

\end{exam}

\pf By substituting  $e^{iw}, e^{ix}, e^{iy}, e^{iz}$ with $w, x, y, z$, respectively,  and then replacing $q^2$ by $q$,  the identity \eqref{add5} becomes
\begin{align}\label{idenwxyz}
&[w^2q,x^2q,y^2q,z^2q;q]_\infty+[-w^2q,-x^2q,-y^2q,-z^2q;q]_\infty\\[6pt]
&\quad =[q xyz/w,q wyz/x,q wxz/y,q wxy/z;q]_\infty
+[-q xyz/w,-q wyz/x,-q wxz/y,-q wxy/z;q]_\infty.\nonumber
\end{align}
Let
\begin{align*}
&\theta_1(w,x,y,z)=[w^2q,x^2q,y^2q,z^2q;q]_\infty, \\[5pt]
&\theta_2(w,x,y,z)=[-w^2q,-x^2q,-y^2q,-z^2q;q]_\infty,\\[5pt]
& \theta_3(w,x,y,z)=[q xyz/w,q wyz/x,q wxz/y,q wxy/z;q]_\infty, \\[5pt]
&\theta_4(w,x,y,z)=[-q xyz/w,-q wyz/x,-q wxz/y,-q wxy/z;q]_\infty.
\end{align*}
Then \eqref{idenwxyz} can be written as
\begin{equation}\label{idenwxyz4}
\theta_4(w,x,y,z)=\theta_1(w,x,y,z)+\theta_2(w,x,y,z)-\theta_3(w,x,y,z).
\end{equation}
Invoking  Theorem \ref{thmexv} and  Theorem  \ref{corexv}, we obtain that for $1\leq k \leq 4$,
\begin{align*}
& \frac{\theta_k(wq^{\frac{1}{2}},xq^{\frac{1}{2}},y,z)}{\theta_k(w,x,y,z)}=\frac{1}{w^2x^2q^2},\\[5pt]
&\frac{\theta_k(wq^{\frac{1}{2}},x,yq^{\frac{1}{2}},z)}{\theta_k(w,x,y,z)}=\frac{1}{w^2y^2q^2},\\[5pt]
&\frac{\theta_k(wq^{\frac{1}{2}},x,y,zq^{\frac{1}{2}})}{\theta_k(w,x,y,z)}=\frac{1}{w^2z^2q^2},\\[5pt]
&\frac{\theta_k(w,xq^{\frac{1}{2}},yq^{\frac{1}{2}},z)}{\theta_k(w,x,y,z)}=\frac{1}{x^2y^2q^2}.
\end{align*}
The set of  exponent vectors of the above contiguous relations is
 \[
 W=\{(2,2,0,0), (2,0,2,0), (2,0,0,2), (0,2,2,0)\}.
 \]
 The fundamental parallelepiped of $W$ is given by
\begin{align*}
\Pi_W &=\{\lambda_1(2,2,0,0)+\lambda_2 (2,0,2,0)+\lambda_3(2,0,0,2)
+\lambda_4(0,2,2,0) \\[5pt]
&\qquad\quad|\, 0\leq\lambda_i <1, i=1,2,3,4\}\cap \mathbb{Z}^4\\[3pt]
&=\big\{(0,0,0,0),(1,1,0,0),(1,0,1,0),(1,0,0,1),(0,1,1,0),(2,1,1,0),\\[3pt]
&\qquad (2,1,0,1),(1,2,1,0),(2,0,1,1),(1,1,2,0),(1,1,1,1),(3,1,1,1),\\[3pt]
&\qquad (2,2,1,1),(2,2,2,0),(2,1,2,1),(3,2,2,1),(1,1,1,0),(2,2,1,0),\\[3pt]
&\qquad (2,1,2,0),(2,1,1,1),(1,2,2,0),(3,2,2,0),(3,2,1,1),(2,3,2,0),\\[3pt]
&\qquad (3,1,2,1),(2,2,3,0),(2,2,2,1),(4,2,2,1),(3,3,2,1),(3,3,3,0),\\[3pt]
&\qquad (3,2,3,1), (4,3,3,1) \big\}.
\end{align*}
Since any nonzero term  of identity \eqref{idenwxyz} has even powers in $w,x,y$ and $z$, \eqref{idenwxyz4} can be further reduced to
\[
\Pi_W'=\{(0,0,0,0), (2,2,2,0)\}.
\]
By Jacobi's triple product identity, we obtain that
\begin{align*}
&[w^0x^0y^0z^0]\,\theta_1=1/(q;q)_\infty^4,\
&& [w^2x^2y^2z^0]\,\theta_1=-q^3/(q;q)_\infty^4,\\[4pt]
&[w^0x^0y^0z^0]\,\theta_2=1/(q;q)_\infty^4, &&[w^2x^2y^2z^0]\,\theta_2=q^3/(q;q)_\infty^4,\\[4pt]
&[w^0x^0y^0z^0]\,\theta_3=1/(q;q)_\infty^4,\ &&[w^2x^2y^2z^0]\,\theta_3=0,\\[4pt]
&[w^0x^0y^0z^0]\,\theta_4=1/(q;q)_\infty^4, &&[w^2x^2y^2z^0]\theta_4=0.
\end{align*}
This completes the proof.  \qed

Let us turn to another addition formula on Jacobi theta functions (see  Whittaker and Watson \cite[p.\,468]{whiwat27}).

\begin{exam} \label{theomt} We have
\begin{align}\label{add1}
2\vartheta_3(w)\vartheta_3(x)\vartheta_3(y)\vartheta_3(z)=&-\vartheta_1(w')\vartheta_1(x')\vartheta_1(y')\vartheta_1(z')
+\vartheta_2(w')\vartheta_2(x')\vartheta_2(y')\vartheta_2(z')\nonumber\\[5pt]
&\quad+\vartheta_3(w')\vartheta_3(x')\vartheta_3(y')\vartheta_3(z')
+\vartheta_4(w')\vartheta_4(x')\vartheta_4(y')
\vartheta_4(z'),
\end{align}
where $w', x', y', z'$ are given by \eqref{wxyz}.
\end{exam}

\pf Substituting  $e^{iw}, e^{ix}, e^{iy}, e^{iz}$ with $w, x, y, z$, respectively,  we may
rewrite  \eqref{add1} as follows:
\begin{align}\label{idenwxyz3}
&2[-qw^2,-qx^2,-qy^2,-qz^2;q^2]_\infty\nonumber\\[4pt]
&\quad=-qwxyz[q^2xyz/w,q^2wyz/x,q^2wxz/y,q^2wxy/z;q^2]_\infty\nonumber\\[4pt]
&\qquad\quad +qwxyz[-q^2xyz/w,-q^2wyz/x,-q^2wxz/y,-q^2wxy/z;q^2]_\infty\nonumber\\[4pt]
&\qquad\quad +[-qxyz/w,-qwyz/x,-qwxz/y,-qwxy/z;q^2]_\infty\nonumber\\[4pt]
&\qquad\quad +[qxyz/w,qwyz/x,qwxz/y,qwxy/z;q^2]_\infty.
\end{align}
Let
\begin{align*}
&\theta_1(w,x,y,z)=wxyz[q^2xyz/w, q^2wyz/x, q^2wxz/y, q^2wxy/z;q^2]_\infty,\\[5pt]
& \theta_2(w,x,y,z)=wxyz[-q^2xyz/w,-q^2wyz/x,-q^2wxz/y,-q^2wxy/z;q^2]_\infty,\\[5pt]
&\theta_3(w,x,y,z)= [-qxyz/w,-qwyz/x,-qwxz/y,-qwxy/z;q^2]_\infty,\\[5pt]
&\theta_4(w,x,y,z)=[qxyz/w,qwyz/x,qwxz/y,qwxy/z;q^2]_\infty,\\[5pt]
&\theta_5(w,x,y,z)=[-qw^2,-qx^2,-qy^2,-qz^2;q^2]_\infty,
\end{align*}
then \eqref{idenwxyz3} becomes
\begin{equation}\label{idenwxyz5}
\theta_5(w,x,y,z)=\frac{1}{2}\big(-q\theta_1(w,x,y,z)+q\theta_2(w,x,y,z)
+\theta_3(w,x,y,z)+\theta_4(w,x,y,z)\big).
\end{equation}
Applying the procedures given in the proofs of Theorem \ref{thmexv} and Theorem \ref{corexv},
 it can be checked that for $1\leq k \leq 5$,
\begin{align*}
& \frac{\theta_k(wq,xq,y,z)}{\theta_k(w,x,y,z)}=\frac{1}{w^2x^2q^2},\\[5pt]
&\frac{\theta_k(wq,x,yq,z)}{\theta_k(w,x,y,z)}=\frac{1}{w^2y^2q^2},\\[5pt]
&\frac{\theta_k(wq,x,y,zq)}{\theta_k(w,x,y,z)}=\frac{1}{w^2z^2q^2},\\[5pt]
&\frac{\theta_k(w,x,yq,zq)}{\theta_k(w,x,y,z)}=\frac{1}{y^2z^2q^2}.
\end{align*}
The set of  exponent vectors of the above contiguous relations is
\[
W=\{(2,2,0,0), (2,0,2,0),(2,0,0,2),(0,0,2,2)\}
\]
and the fundamental parallelepiped of $W$ equals
\begin{align*}
\Pi_W &=\{\lambda_1(2,2,0,0)+\lambda_2 (2,0,2,0)+\lambda_3(2,0,0,2)
+\lambda_4(0,0,2,2) \ |\ 0\leq\lambda_i <1, 1\leq i \leq 4\}\cap \mathbb{Z}^4\\[3pt]
&=\big\{ (0,0,0,0),(1,1,0,0),(1,0,1,0),(1,0,0,1),(0,0,1,1),(2,1,1,0),
(2,1,0,1),\\[3pt]
&\qquad (1,1,1,1),(2,0,1,1),(1,0,2,1),(1,0,1,2),(3,1,1,1),
(2,1,1,2),(2,1,2,1),\\[3pt]
&\qquad(2,0,2,2),(3,1,2,2),(1,0,1,1), (2,1,1,1),
(2,0,2,1),(2,0,1,2),(1,0,2,2),\\[3pt]
&\qquad(3,1,2,1),(3,1,1,2),(2,1,2,2),
(3,0,2,2),(2,0,3,2), (2,0,2,3),(4,1,2,2),\\[3pt]
&\qquad(3,1,2,3),(3,1,3,2), (3,0,3,3),(4,1,3,3) \big\}.
\end{align*}
Noting that  \eqref{idenwxyz3} contains terms only  with even powers in $w,x,y,z$,  \eqref{idenwxyz5} can be further reduced to  equalities on the coefficients with exponent vectors in
\[
\Pi_W'=\{(0,0,0,0),(2,0,2,2)\}.
\]
These two relations are just consequences of Jacobi's triple product identity,
and so the proof is complete.   \qed

To conclude this paper, we  give an  example showing that
our approach does not restrict  to identities on multiple theta functions  in the form of \eqref{thetaeq}.
Let us consider the following identity:
\begin{align}\label{wwyz}
\vartheta_1(y+z)\vartheta_1(y-z)\vartheta_4^2=\vartheta_3^2(y)\vartheta_2^2(z)-
\vartheta_2^2(y)\vartheta_3^2(z),
\end{align}
where $\vartheta_4=\vartheta_4(0)$  (see Whittaker and Watson \cite{whiwat27}).

{\noindent \it Proof.}
Substituting $e^{iy}$ and  $e^{iz}$ with $\sqrt{y}$ and
$\sqrt{z}$, respectively, \eqref{wwyz} can be rewritten as
\begin{align}\label{wwsquare}
&(q;q^2)_\infty^4 y[q^2yz,q^2y/z;q^2]_\infty
=-z[-qy,-q^2z;q^2]_\infty^2+y
[-q^2y,-qz;q^2]_\infty^2.
\end{align}
Let
\begin{align*}
&\theta_1(y,z)=z[-qy,-q^2z;q^2]_\infty^2,\\[5pt]
&\theta_2(y,z)=y[-q^2y,-qz;q^2]_\infty^2,\\[5pt]
&\theta_3(y,z)=y[q^2yz,q^2y/z;q^2]_\infty.
\end{align*}
Then (\ref{wwsquare}) can be written as
\[
\theta_3(y,z)=\frac{1}{(q;q^2)_\infty^4}\big(-\theta_1(y,z)+\theta_2(y,z)\big).
\]
Even though our approach do not apply to the multiple
theta functions $\theta_1$ and $\theta_2$,
it is still possible to find two contiguous relations with linearly independent
exponent vectors satisfied by
$\theta_1$, $\theta_2$, and $\theta_3$. It can be easily checked that for $1\leq k \leq 3$,
\begin{align}\label{contigwwsquare}
& \frac{\theta_k(yq^2,z)}{\theta_k(y,z)}=\frac{1}{y^2q^2},\nonumber \\[5pt]
&\frac{\theta_k(y,zq^2)}{\theta_k(y,z)}=\frac{1}{z^2q^2}.
\end{align}
Therefore, Stanley's Lemma can be used to reduce  the coefficients of $\theta_k$ to a finite number of initial coefficients based on the contiguous relations in \eqref{contigwwsquare}.
The set of  exponent vectors of the two contiguous relations is
\[
W=\{(2,0), (0,2)\}
\]
and the fundamental parallelepiped of $W$ is given by
\begin{align*}
\Pi_W &=\{\lambda_1(2,0)+\lambda_2 (0,2) \ |\ 0\leq\lambda_1,\lambda_2
<1\}\cap \mathbb{Z}^2\\[6pt]
&=\{(0,0),(1,0), (0,1),(1,1)\}.
\end{align*}
Because of the symmetry of $y$ and $z$, $\Pi_W$ can be further reduced to
\[
\Pi'_W=\{(0,0), (1,0), (1,1)\}.
\]
By Jacobi's triple product identity, we obtain that
\begin{align*}
&[y^0z^0]\theta_1=\frac{2(-q^2;q^2)_\infty^2(q^4;q^4)_\infty^2}{(q^2;q^2)_\infty^4},
&& [y^0z^0]\theta_2=\frac{2(-q^2;q^2)_\infty^2(q^4;q^4)_\infty^2}{(q^2;q^2)_\infty^4},\\[5pt]
&[y^0z^0]\theta_3=0, &&[y^1z^0]\theta_1=\frac{4q(-q^4;q^4)_\infty^4(q^4;q^4)_\infty^2}{(q^2;q^2)_\infty^4},\\[5pt]
&
[y^1z^0]\theta_2=\frac{(-q^2;q^2)_\infty^2(q^4;q^4)_\infty^2}{(q^2;q^2)_\infty^4},
&& [y^1z^0]\theta_3=\frac{1}{(q^2;q^2)_\infty^2},\\[5pt]
&[y^1z^1]\theta_1=\frac{2q(-q^2;q^2)_\infty^2 (q^4;q^4)_\infty^2}{(q^2;q^2)_\infty^4},
&&
[y^1z^1]\theta_2=\frac{2q(-q^2;q^2)_\infty^2(q^4;q^4)_\infty^2}{(q^2;q^2)_\infty^4}, \\[5pt]
& [y^1z^1]\theta_3=0.
\end{align*}
The equality of the coefficients of $y^1z^0$ simplifies to
Ewell's identity  \cite[Eq.\,(2.3)]{ewell98}
\begin{equation*}
(-q^2;q^4)_\infty^4-4q(-q^4;q^4)_\infty^4=(q^2;q^4)_\infty^2(q;q^2)_\infty^4.
\end{equation*} Thus the proof of \eqref{wwyz} is complete. \qed

\vskip 15pt \noindent {\small {\bf Acknowledgments.}  We wish to
thank  Doron Zeilberger, Yi Zhang and the referee for valuable suggestions. This work was supported by the 973 Project,  the National Science Foundation, and the Natural Science Foundation of Tianjin City, China.

\end{document}